# Distribution-Free Tests of Independence in High Dimensions

BY FANG HAN

*Department of Statistics, University of Washington, Box 354322, Seattle, Washington 98195, U.S.A.*

fanghan@uw.edu

SHIZHE CHEN

*Department of Statistics, Columbia University, Room 1005 SSW, MC 4690, 1255 Amsterdam Avenue, New York, New York 10027, U.S.A.*

shizhe.chen@gmail.com

AND HAN LIU

*Department of Operations Research and Financial Engineering, Princeton University, Sherrerd Hall, Charlton Street, Princeton, New Jersy 08544, U.S.A.*

hanliu@princeton.edu

SUMMARY

We consider the testing of mutual independence among all entries in a $d$-dimensional random vector based on $n$ independent observations. We study two families of distribution-free test statistics, which include Kendall's tau and Spearman's rho as important examples. We show that under the null hypothesis the test statistics of these two families converge weakly to Gumbel distributions, and propose tests that control the type I error in the high-dimensional setting where $d > n$. We further show that the two tests are rate-optimal in terms of power against sparse alternatives, and outperform competitors in simulations, especially when $d$ is large.

*Some key words*: Gumbel distribution; Kendall's tau; Linear rank statistic; Mutual independence; Rank-type $U$-statistic; Spearman's rho.

## 1. INTRODUCTION
### 1·1. *Literature review*

Consider a $d$-dimensional continuous random vector $X$, where $X = (X_1, \ldots, X_d)^{\mathrm{T}} \in \mathcal{R}^d$. Given $n$ samples, we aim to test the null hypothesis $H_0 : X_1, \ldots, X_d$ are mutually independent.

This problem has been intensively studied when $X$ is multivariate Gaussian. When $d < n$, methods proposed include the likelihood ratio test (Anderson, 2003), Roy's (1957) largest root test, and Nagao's (1973) test, which test the identity of the Pearson's covariance matrix $\Sigma$ or correlation matrix $R$ using their sample counterparts. When $d$ and $n$ both grow and the ratio $d/n$ does not converge to zero, classic likelihood ratio tests perform poorly since the sample eigenvalues do not converge to their population counterparts (Bai & Yin, 1993). This has motivated work in high-dimensional settings.

In what follows, let $\gamma$ denote the limit of $d/n$ as $n$ and $d$ diverge to infinity. When $0 < \gamma \leq 1$, Bai et al. (2009) and Jiang & Yang (2013) propose, and establish asymptotic normality of, cor-rected likelihood ratio test statistics. Specifically, Bai et al. (2009) consider the regime $\gamma \in (0, 1)$, and Jiang et al. (2012) extend it to the case $\gamma = 1$. Johnstone (2001) and Bao et al. (2012) prove the Tracy–Widom law for the null limiting distributions of Roy's largest root test statistics. Bao et al. (2012)'s result is only valid when $\gamma \in (0, 1)$, while the result in Johnstone (2001) applies to the case $\gamma = 1$. These results are further generalized to $\gamma > 1$ in Péché (2009) and Pillai & Yin (2012), with possibly non-Gaussian observations. When $\gamma$ can be arbitrarily large but bounded, Ledoit & Wolf (2002) and Schott (2005) propose, and establish asymptotic normality of, corrected Nagao test statistics. Jiang (2004) proposes a test statistic based on the largest sample correlation coefficient and shows that it converges to a Gumbel, also known as an extreme-value type I, distribution. With some adjustments, Birke & Dette (2005) and Cai & Jiang (2012) prove that the tests of Ledoit & Wolf (2002) and Jiang (2004) extend to the case $\gamma = \infty$. To the best of our knowledge, there is no result generalizing the test of Schott (2005) to the regime $\gamma = \infty$.

When $\gamma$ can equal infinity, Srivastava (2006) proposes a corrected likelihood ratio test using only nonzero sample eigenvalues. Srivastava (2005) introduces a test using unbiased estimators of the traces of powers of the covariance matrix. Cai & Ma (2013) show that the test of Chen et al. (2010) uniformly dominates the corrected likelihood ratio tests of Bai et al. (2009) and Jiang & Yang (2013); the three test statistics are asymptotically normal. Zhou (2007) modifies Jiang (2004)'s test and shows that the null limiting distribution of the test statistic is Gumbel.

Most of the aforementioned tests are valid only under normality. For non-Gaussian data, testing $H_0$ in high dimensions is not as well studied: Péché (2009) and Pillai & Yin (2012) study Roy's largest root test for sub-Gaussian data; Bao et al. (2015) study the Spearman's rho statistic; and Jiang (2004) studies the largest off-diagonal entry in the sample correlation matrix. In particular, Jiang (2004) shows that, for testing a simplified version of $H_0$, the normality assumption can be relaxed to $E(|X|^r) < \infty$ for some $r > 30$. Later, Zhou (2007) modifies Jiang (2004)'s test to require only $r \geq 6$. See also Li & Rosalsky (2006), Zhou (2007), Liu et al. (2008), Li et al. (2010), Cai & Jiang (2011), Cai & Jiang (2012), Shao & Zhou (2014), and Han et al. (2017).

This paper investigates testing $H_0$ in high dimensions. The asymptotic regime of interest is when $d$ and $n$ both grow and $d/n$ can diverge or converge to any non-negative value. Our main focus is on non-parametric rank-based tests and their optimality. We consider two families of rank-based test statistics including Spearman's rho (Spearman, 1904) and Kendall's tau (Kendall, 1938), and prove that under the null hypothesis they converge weakly to Gumbel distributions. We also provide power analysis and establish optimality of the proposed tests against sparse alternatives, explicitly defined in Section 4. In particular, we show that the tests based on Spearman's rho and Kendall's tau are rate-optimal against sparse alternatives.

### 1·2. *Other related work*

Testing $H_0$ is related to testing bivariate independence. To test the independence between two random variables taking scalar values, Hotelling & Pabst (1936) and Kendall (1938) propose to use Spearman's rho and Kendall's tau statistics, and Hoeffding (1948b) proposes the $D$ statistic. To test the independence of two random vectors with possibly very high dimensions, Bakirov et al. (2006), Székely & Rizzo (2013), and Jiang et al. (2013) propose tests based on normalized distance between characteristic functions, distance correlations, and modified likelihood ratios. However, we cannot directly apply these results to test $H_0$ without multiple testing adjustments.

A notable alternative to Pearson's correlation coefficient is Spearman's rho. Zhou (2007) establishes the limiting distribution of the largest off-diagonal entry of the Spearman's rho correlation matrix, but does not provide a power analysis of the corresponding test. This paper includes the result in Zhou (2007) as a special case.

Many have considered testing independence based on kernel methods (Gretton et al., 2007; Fukumizu et al., 2008; Póczos et al., 2012; Reddi & Póczos, 2013). They focus on the kernel dependence measures, the Hilbert–Schmidt norm of the cross-covariance operator (Gretton et al., 2007), or the normalized cross-covariance operator (Fukumizu et al., 2008). Using these dependence measures, early works consider testing independence between two random variables (Gretton et al., 2007; Fukumizu et al., 2008) that might live in arbitrary sample spaces. Recently, Reddi & Póczos (2013) generalize the proposal in Fukumizu et al. (2008) and propose a copula-based kernel dependence measure for testing mutual independence. Póczos et al. (2012) offer an alternative kernel-based test using the maximum mean discrepancy (Borgwardt et al., 2006) between the empirical copula and the joint distribution of $d$ independent uniform random variables. However, existing kernel-based tests are proposed in the low-dimensional setting; Ramdas et al. (2015) have shown that the kernel-based tests have low power in high dimensions.

During the preparation of this paper, it has come to our attention that recent works by Mao (2016) and Leung & Drton (2017) also consider testing $H_0$ for non-Gaussian data in high dimensions. They propose tests based on sums of rank correlations, such as Kendall's tau (Leung & Drton, 2017) and Spearman's rho (Mao, 2016; Leung & Drton, 2017). They further establish asymptotic normality of the proposed test statistics in the case when $\gamma$ can be arbitrarily large but bounded. In particular, the theory in Mao (2016) follows from the procedure developed in Schott (2005), and the theory in Leung & Drton (2017) relies on $U$-statistics theory.

## 2. TESTING PROCEDURES

### 2·1. *Two families of tests*

Let $\{X_{i,\cdot} : X_{i,\cdot} = (X_{i,1}, \ldots, X_{i,d})^{\mathrm{T}}, i = 1, \ldots, n\}$ be $n$ independent replicates of a $d$-dimensional random vector $X \in \mathcal{R}^d$. To avoid discussion of possible ties, we consider continuous random vectors. For any two entries $j \neq k \in \{1, \ldots, d\}$, let $Q_{ni}^j$ be the rank of $X_{i,j}$ in $\{X_{1,j}, \ldots, X_{n,j}\}$ and let $R_{ni}^{jk}$ be the relative rank of the $k$-th entry corresponding to the $j$-th entry; that is, $R_{ni}^{jk} \equiv Q_{ni'}^k$ subject to the constraint that $Q_{ni'}^j = i$, for $i = 1, \ldots, n$.

We propose two families of non-parametric tests based on the relative ranks. The first family includes tests based on simple linear rank statistics of the form

$$V_{jk} \equiv \sum_{i=1}^{n} c_{ni} g\{R_{ni}^{jk}/(n+1)\} \quad (j \neq k \in \{1, \ldots, d\}), \tag{1}$$

where $\{c_{ni}, i = 1, \ldots, n\}$ form an array of constants called the regression constants, and $g(\cdot)$ is a Lipschitz function called the score function. We assume $\sum_{i=1}^{n} c_{ni}^2 > 0$ to avoid triviality. It is immediately clear that Spearman's rho belongs to the family of simple linear rank statistics. For accommodating tests of independence, we further pose the alignment assumption,

$$c_{ni} = n^{-1} f\{i/(n+1)\}, \tag{2}$$

where $f(\cdot)$ is a Lipschitz function. Under this assumption, the simple linear rank statistic is a general measure of the agreements between the ranks of two sequences. It will be made clear, in Sections 3 and 4, that the alignment assumption (2) is not required in deriving the null limiting distribution, but is crucial in the power analysis.

The second family includes tests based on rank-type $U$-statistics, which are $U$-statistics of order $m < n$ that depend only on relative ranks $\{R_{ni}^{jk}, i = 1, \ldots, n\}$. In other words, a rank-type

$U$-statistic takes the form

$$U_{jk} \equiv \frac{1}{n(n-1)\cdots(n-m+1)} \sum_{i_1 \neq i_2 \neq \cdots \neq i_m} h(X_{i_1,\{j,k\}}, \ldots, X_{i_m,\{j,k\}}) \quad (j \neq k \in \{1,\ldots,d\}) \quad (3)$$

and $U_{jk}$ only depends on $\{R_{ni}^{jk}\}_{i=1}^n$. Here for any vector $X_{i,\cdot}$ and some index set $\mathcal{A} \subset \{1,\ldots,d\}$, we let $X_{i,\mathcal{A}}$ be the sub-vector of $X_{i,\cdot}$ with entries in the index set $\mathcal{A}$. The kernel function $h(\cdot)$ is assumed to be bounded, but not necessarily symmetric. The boundedness assumption is mild since correlation is the object of interest.

Next, we propose two tests based on the two families of statistics, respectively. We begin with a testing procedure based on simple linear rank statistics. Under $H_0$, the distribution of $V_{jk}$ is irrelevant to the specific distribution of $X$ for all $j \neq k \in \{1,\ldots,d\}$. Accordingly, the mean and variance of $V_{jk}$ are calculable without knowing the true distribution. Let $E_{H_0}(\cdot)$ and $\text{var}_{H_0}(\cdot)$ be the expectation and variance of a certain statistic under $H_0$. We have

$$E_{H_0}(V_{jk}) = \bar{g}_n \sum_{i=1}^n c_{ni}, \quad \text{var}_{H_0}(V_{jk}) = \frac{1}{n-1} \sum_{i=1}^n \left[g\{i/(n+1)\} - \bar{g}_n\right]^2 \sum_{i=1}^n (c_{ni} - \bar{c}_n)^2, \quad (4)$$

where $\bar{g}_n \equiv n^{-1} \sum_{i=1}^n g\{i/(n+1)\}$ is the sample mean of $g\{R_{ni}^{jk}/(n+1)\} (i=1,\ldots,n)$. Based on $\{V_{jk}, 1 \leq j < k \leq d\}$, we propose the following statistic for testing $H_0$:

$$L_n \equiv \max_{j<k} |V_{jk} - E_{H_0}(V_{jk})|. \quad (5)$$

As with simple linear rank statistics, the expectation and variance of the rank-type $U$-statistics, $E_{H_0}(U_{jk})$ and $\text{var}_{H_0}(U_{jk})$, can be calculated analytically. We can test $H_0$ using

$$\widetilde{L}_n \equiv \max_{j<k} |U_{jk} - E_{H_0}(U_{jk})|. \quad (6)$$

Detailed studies of $L_n$ and $\widetilde{L}_n$'s null limiting distributions are deferred to Section 3. Instead, we give some intuition here. Under certain conditions, the standardized version of $V_{jk}$ or $U_{jk}$ is asymptotically normal. Accordingly, the standardized version of $L_n^2$ or $\widetilde{L}_n^2$, is asymptotically close to the maximum of $d(d-1)/2$ independent chi-squared random variables with degree of freedom one. The latter converges weakly to a Gumbel distribution after adjustment.

Let $\sigma_V^2$ and $\sigma_U^2$ be the variances of $n^{1/2} V_{jk}$ and $n^{1/2} U_{jk}$ under $H_0$:

$$\sigma_V^2 \equiv n\text{var}_{H_0}(V_{jk}), \quad \sigma_U^2 \equiv n\text{var}_{H_0}(U_{jk}). \quad (7)$$

We propose the size-$\alpha$ tests $T_\alpha$ and $\widetilde{T}_\alpha$ of $H_0$ as follows:

$$T_\alpha \equiv I\left(\frac{nL_n^2}{\sigma_V^2} - 4\log d + \log\log d \geq q_\alpha\right), \widetilde{T}_\alpha \equiv I\left(\frac{n\widetilde{L}_n^2}{\sigma_U^2} - 4\log d + \log\log d \geq q_\alpha\right). \quad (8)$$

Here $I(\cdot)$ represents the indicator function and

$$q_\alpha \equiv -\log(8\pi) - 2\log\log(1-\alpha)^{-1} \quad (9)$$

is the $1-\alpha$ quantile of the Gumbel distribution function $\exp\{-(8\pi)^{-1/2} \exp(-y/2)\}$. In the sequel, we only consider a fixed nominal significance level, e.g., $\alpha = 0 \cdot 05$.

As an alternative, we can simulate the exact distribution of the studied statistic and choose $q_\alpha$ to be the $1-\alpha$ quantile of the corresponding empirical distribution. The simulation-based approach is discussed in the Supplementary Material.

2·2. *Examples*

In the following, we provide four distribution-free tests of independence that belong to the two general families defined in Section 2·1.

*Example* 1 (*Spearman's rho*). Recall that $Q_{ni}^j$ and $Q_{ni}^k$ are the ranks of $X_{i,j}$ and $X_{i,k}$ among $\{X_{1,j},\ldots,X_{n,j}\}$ and $\{X_{1,k},\ldots,X_{n,k}\}$, respectively. Spearman's rho is defined as

$$\rho_{jk} = \frac{\sum_{i=1}^n (Q_{ni}^j - \bar{Q}_n^j)(Q_{ni}^k - \bar{Q}_n^k)}{\left\{\sum_{i=1}^n (Q_{ni}^j - \bar{Q}_n^j)^2 \sum_{i=1}^n (Q_{ni}^k - \bar{Q}_n^k)^2\right\}^{1/2}}$$
$$= \frac{12}{n(n^2-1)} \sum_{i=1}^n \left(i - \frac{n+1}{2}\right)\left(R_{ni}^{jk} - \frac{n+1}{2}\right) \quad (j \neq k \in \{1,\ldots,d\}), \qquad (10)$$

where $\bar{Q}_n^j = \bar{Q}_n^k \equiv (n+1)/2$. This is a simple linear rank statistic, and

$$E_{H_0}(\rho_{jk}) = 0, \quad \text{var}_{H_0}(\rho_{jk}) = (n-1)^{-1} \quad (j \neq k \in \{1,\ldots,d\}).$$

According to (8), the corresponding test statistic is

$$T_\alpha^\rho = I\{(n-1)\max_{j<k} \rho_{jk}^2 - 4\log d + \log\log d \geq q_\alpha\}.$$

*Example* 2 (*Kendall's tau*). Kendall's tau is defined as, for $j \neq k \in \{1,\ldots,d\}$,

$$\tau_{jk} = \frac{2}{n(n-1)} \sum_{i<i'} \text{sign}(X_{i',j} - X_{i,j})\text{sign}(X_{i',k} - X_{i,k}) = \frac{2}{n(n-1)} \sum_{i<i'} \text{sign}(R_{ni'}^{jk} - R_{ni}^{jk}),$$

where the sign function $\text{sign}(\cdot)$ is defined as $\text{sign}(x) = x/|x|$, with the convention $0/0 = 0$. This statistic is a function of the relative ranks $\{R_{ni}^{jk}, i = 1,\ldots,n\}$ and also a $U$-statistic with a bounded kernel $h(x_{1,\{1,2\}}, x_{2,\{1,2\}}) \equiv \text{sign}(x_{1,1} - x_{2,1})\text{sign}(x_{1,2} - x_{2,2})$. Accordingly, Kendall's tau is a rank-type $U$-statistic. Moreover,

$$E_{H_0}(\tau_{jk}) = 0, \quad \text{var}_{H_0}(\tau_{jk}) = \frac{2(2n+5)}{9n(n-1)} \quad (j \neq k \in \{1,\ldots,d\}).$$

According to (8), the proposed test statistic based on Kendall's tau is

$$T_\alpha^\tau = I\left\{\frac{9n(n-1)}{2(2n+5)} \max_{j<k} \tau_{jk}^2 - 4\log d + \log\log d \geq q_\alpha\right\}.$$

*Example* 3 (*A major part of Spearman's rho*). Although Spearman's rho is not a $U$-statistic, by Hoeffding (1948a), we can write, for $j \neq k \in \{1,\ldots,d\}$,

$$\rho_{jk} = \frac{n-2}{n+1}\widehat{\rho}_{jk} + \frac{3\tau_{jk}}{n+1}, \qquad (11)$$

where

$$\widehat{\rho}_{jk} = \frac{3}{n(n-1)(n-2)} \sum_{i \neq i' \neq i''} \text{sign}(X_{i,j} - X_{i',j})\text{sign}(X_{i,k} - X_{i'',k}).$$

Here $\widehat{\rho}_{jk}$ is a $U$-statistic with degree three and an asymmetric bounded kernel. Moreover,

$$E_{H_0}(\widehat{\rho}_{jk}) = 0, \quad \text{var}_{H_0}(\widehat{\rho}_{jk}) = \frac{n^2 - 3}{n(n-1)(n-2)} \quad (j \neq k \in \{1,\ldots,d\}).$$

As in (8), we propose the test based on $\{\widehat{\rho}_{jk}, 1 \leq j < k \leq d\}$ as

$$T_\alpha^{\widehat{\rho}} = I\left\{\frac{n(n-1)(n-2)}{n^2-3}\max_{j<k}\widehat{\rho}_{jk}^2 - 4\log d + \log\log d \geq q_\alpha\right\}.$$

*Example* 4 (*Projection of Kendall's tau to simple linear rank statistics*). Kendall's tau does not belong to the family of simple linear rank statistics. However, by the projection argument in Hájek (1968), $\tau_{jk}$ can be approximated by

$$\widehat{\tau}_{jk} = \frac{8}{n^2(n-1)}\sum_{i=1}^{n}\left(i - \frac{n+1}{2}\right)\left(R_{ni}^{jk} - \frac{n+1}{2}\right) \quad (j \neq k \in \{1,\ldots,d\}).$$

Using the variance of $\rho_{jk}$ and the relationship between $\rho_{jk}$ and $\widehat{\tau}_{jk}$, it is easy to obtain

$$E_{H_0}(\widehat{\tau}_{jk}) = 0, \quad \mathrm{var}_{H_0}(\widehat{\tau}_{jk}) = \frac{4(n+1)^2}{9n^2(n-1)} \quad (j \neq k \in \{1,\ldots,d\}).$$

We observe that $\mathrm{var}_{H_0}(\widehat{\tau}_{jk})/\mathrm{var}_{H_0}(\tau_{jk})$ goes to unity as $n$ grows, indicating that $\widehat{\tau}_{jk}$ is asymptotically equivalent to $\tau_{jk}$ under $H_0$. The proposed test statistic is

$$T_\alpha^{\widehat{\tau}} = I\left\{\frac{9n^2(n-1)}{4(n+1)^2}\max_{j<k}\widehat{\tau}_{jk}^2 - 4\log d + \log\log d \geq q_\alpha\right\}.$$

*Remark* 1. We have considered two families of test statistics: that of simple linear rank statistics and that of rank-type $U$-statistics. Waerden (1957) and Woodworth (1970) studied the performance of Spearman's rho and Kendall's tau in testing bivariate independence under normality, and show that Spearman's rho is more efficient than Kendall's tau when $n$ is small, while the reverse is true if $n$ is large. Although the threshold point is theoretically calculable, in practice it is very difficult to approximate it.

## 3. LIMITING NULL DISTRIBUTIONS

This section characterizes the limiting distributions of $L_n$ and $\widetilde{L}_n$ under $H_0$. We first introduce some necessary notation. Let $v = (v_1,\ldots,v_d)^{\mathrm{T}} \in \mathcal{R}^d$ be a $d$-vector and let $M = [M_{jk}] \in \mathcal{R}^{d\times d}$ be a $d \times d$ matrix. For any index sets $I, J \subset \{1,\ldots,d\}$, let $v_I$ be the sub-vector of $v$ with entries indexed by $I$, and $M_{I,J}$ be the sub-matrix of $M$ with rows indexed by $I$ and columns indexed by $J$. Let $\lambda_{\min}(M)$ denote the smallest eigenvalue of $M$. For two sequences $\{a_1, a_2, \ldots\}$ and $\{b_1, b_2, \ldots\}$, we write $a_n = O(b_n)$ if there exists some constant $C$ such that, for any sufficiently large $n$, $|a_n| \leq C|b_n|$. We write $a_n = o(b_n)$, if for any positive constant $c$ and sufficiently large $n$, $|a_n| \leq c|b_n|$. We write $a_n = o_y(b_n)$, if the constant depends on some scalar $y$, i.e., $|a_n| \leq c_y|b_n|$. We study the asymptotics of triangular arrays (Greenshtein & Ritov, 2004), allowing the dimension $d \equiv d_n$ to grow with $n$. We use $c$ and $C$ to represent generic positive constants, whose values may vary at different locations.

We first consider the simple linear rank statistic $V_{jk}$. The following theorem shows that, under $H_0$ and some regularity conditions on the regression constants $\{c_{n1},\ldots,c_{nn}\}$, the statistic $nL_n^2/\sigma_V^2 - 4\log d + \log\log d$ converges weakly to a Gumbel distribution.

THEOREM 1. *Suppose that the simple linear rank statistics $\{V_{jk}, 1 \leq j < k \leq d\}$ take the form* (1) *with regression constants $\{c_{n1}, \ldots, c_{nn}\}$ satisfying*

$$\max_{1\leq i\leq n} |c_{ni}-\bar{c}_n|^2 \leq \frac{C_1^2}{n}\sum_{i=1}^n (c_{ni}-\bar{c}_n)^2, \quad \left|\sum_{i=1}^n (c_{ni}-\bar{c}_n)^3\right|^2 \leq \frac{C_2^2}{n}\left\{\sum_{i=1}^n (c_{ni}-\bar{c}_n)^2\right\}^3, \quad (12)$$

*where $\bar{c}_n \equiv \sum_{i=1}^n c_{ni}$ represents the sample mean of the regression constants and $C_1, C_2$ are two constants. Further suppose that the score function $g(\cdot)$ is differentiable with bounded Lipschitz constant. We then have, under $H_0$, if $\log d = o(n^{1/3})$ as $n$ grows, then for any $y \in \mathcal{R}$,*

$$|\mathrm{pr}(nL_n^2/\sigma_V^2 - 4\log d + \log\log d \leq y) - \exp\{-(8\pi)^{-1/2}\exp(-y/2)\}| = o_y(1),$$

*where $L_n$ and $\sigma_V^2$ are defined in* (5) *and* (7).

In Theorem 1, conditions in the form (12) are common for the simple linear rank statistics to be asymptotically normal or to deviate moderately from normality; see Hájek et al. (1999) and Kallenberg (1982). Seoh et al. (1985) propose similar conditions for $\{c_{ni}, i = 1, \ldots, n\}$. The Lipschitz condition rules out the Fisher–Yates statistic, where $g(\cdot)$ is proportional to $\Phi^{-1}\{\cdot/(n+1)\}$ and $\Phi^{-1}(\cdot)$ represents the quantile function of the standard Gaussian.

Theorem 1 gives a distribution-free result for testing $H_0$ (see Chapter 31 in Kendall & Stuart, 1961). In contrast, tests based on sample covariance and correlation matrices (Jiang, 2004; Li et al., 2010; Cai & Jiang, 2011; Shao & Zhou, 2014) are not distribution-free: for instance, Li et al. (2010) and Shao & Zhou (2014) impose moment requirements on $X$.

Spearman's rho is a simple linear rank statistic, and satisfies the conditions in (12). Therefore, Theorem 1 is a strict generalization of Theorem 1.2 in Zhou (2007).

We then turn to rank-type $U$-statistics. The next theorem mirrors Theorem 1.

THEOREM 2. *Suppose that the rank-type $U$-statistics $\{U_{jk}, 1 \leq j < k \leq d\}$ are of the form* (3), *of degree $m$, and the kernel function $h(\cdot)$ is bounded and non-degenerate. We then have, under $H_0$, if $\log d = o(n^{1/3})$ as $n$ grows, then for any $y \in \mathcal{R}$,*

$$\left|\mathrm{pr}(n\widetilde{L}_n^2/\sigma_U^2 - 4\log d + \log\log d \leq y) - \exp\{-(8\pi)^{-1/2}\exp(-y/2)\}\right| = o_y(1),$$

*where $\widetilde{L}_n$ and $\sigma_U^2$ are defined in* (6) *and* (7).

The assumption on $h(\cdot)$ states that the rank-type $U$-statistic is non-degenerate, and hence rules out Hoeffding's $D$ statistic.

Corollary 1 shows that the tests $T_\alpha$ and $\widetilde{T}_\alpha$ can effectively control the size.

COROLLARY 1. *Suppose that the conditions in Theorems* 1 *or* 2 *hold, respectively, then*

$$\mathrm{pr}(T_\alpha = 1 \mid H_0) = \alpha + o(1), \quad \mathrm{pr}(\widetilde{T}_\alpha = 1 \mid H_0) = \alpha + o(1).$$

Furthermore, all test statistics in Examples 1–4 converge weakly to a Gumbel distribution.

COROLLARY 2. *Under the regime $\log d = o(n^{1/3})$ as $n$ grows,*

$$\mathrm{pr}(T_\alpha^\ell = 1 \mid H_0) = \alpha + o(1) \quad (\ell \in \{\rho, \tau, \widehat{\rho}, \widehat{\tau}\}),$$

*where $T_\alpha^\ell$ corresponds to the test statistics introduced in Examples* 1–4 *for $\ell \in \{\rho, \tau, \widehat{\rho}, \widehat{\tau}\}$.*

## 4. POWER ANALYSIS AND OPTIMALITY PROPERTIES

### 4·1. *Sparse alternatives*

Let $\mathcal{U}(c)$ be a set of matrices indexed by a constant $c$

$$\mathcal{U}(c) \equiv \Big\{ M \in \mathcal{R}^{d \times d} : \mathrm{diag}(M) = I_d, M = M^{\mathrm{T}}, \max_{1 \leq j < k \leq d} |M_{jk}| \geq c(\log d/n)^{1/2} \Big\}, \quad (13)$$

where $I_d$ represents the $d \times d$ identity matrix, and $\mathrm{diag}(M)$ represents a matrix with diagonals equal the diagonals of $M$ and all off-diagonals equal zero.

We define the random matrix $\widehat{V} = [\widehat{V}_{jk}] \in \mathcal{R}^{d \times d}$ as

$$\widehat{V}_{jk} = \widehat{V}_{kj} = \sigma_V^{-1} \{ V_{jk} - E_{H_0}(V_{jk}) \}, \quad \widehat{V}_{\ell\ell} = 1 \quad (1 \leq j < k \leq d, 1 \leq \ell \leq d),$$

where $\sigma_V$ is defined in (7) and $\{V_{jk}, 1 \leq j < k \leq d\}$ are the simple linear rank statistics. Let population version of $\widehat{V}$ be $V \equiv E(\widehat{V})$. We study the power of tests against the alternative

$$H_a^V(c) \equiv \{F(X) : V\{F(X)\} \in \mathcal{U}(c)\}, \quad (14)$$

where $F(X)$ is the joint distribution function of $X$ and we write $V\{F(X)\}$ to emphasize that $V = E(\widehat{V}) = \int \widehat{V} \mathrm{d} F(X)$ is a function of $F(X)$.

Similarly, we define the random matrix $\widehat{U} = [\widehat{U}_{jk}] \in \mathcal{R}^{d \times d}$ as

$$\widehat{U}_{jk} = \widehat{U}_{kj} = \frac{U_{jk} - E_{H_0}(U_{jk})}{\widetilde{\sigma}_U}, \quad \widehat{U}_{\ell\ell} = 1 \quad (1 \leq j < k \leq d; 1 \leq \ell \leq d),$$

where $\{U_{jk}, 1 \leq j < k \leq d\}$ are the rank-type $U$-statistics and

$$\widetilde{\sigma}_U^2 \equiv m^2 \mathrm{var}_{H_0} \big[ E_{H_0} \{ h(X_{1,\{1,2\}}, \ldots, X_{m,\{1,2\}}) \mid X_{1,\{1,2\}} \} \big]. \quad (15)$$

We define the population version of $\widehat{U}$ to be $U \equiv E(\widehat{U})$. Then, we study the power of tests against the alternative

$$H_a^U(c) \equiv \{F(X) : U\{F(X)\} \in \mathcal{U}(c)\}. \quad (16)$$

When studying rate-optimality, we consider the following alternative

$$H_a^R(c) \equiv \{F(X) : R\{F(X)\} \in \mathcal{U}(c)\}, \quad (17)$$

where $R$ is the population correlation matrix. Section 4·3 clarifies why we use $H_a^R(c)$ in (17).

All three alternatives are based on the set of matrices $\mathcal{U}(c)$, of which at least one entry's magnitude is larger than $C(\log d/n)^{1/2}$ for some large constant $C$, so we call (14), (16), and (17) the sparse alternatives.

The three alternatives may not be equivalent. For instance, $\max_{1 \leq j < k \leq d} |V_{jk}| \geq c(\log d/n)^{1/2}$ does not imply $\max_{1 \leq j < k \leq d} |U_{jk}| \geq c(\log d/n)^{1/2}$. The exact relationship between $H_a^V$ and $H_a^U$ is intriguing. Taking Kendall's tau and Spearman's rho as examples, Fredricks & Nelsen (2007a) show that, for a bivariate random vector, the ratio between the population analogs of the two statistics converges to $3/2$ as the joint distribution approaches independence. Under a fixed alternative, however, the relationship between the population analogs of Kendall's tau and Spearman's rho remains unclear and probably depends heavily on the specific distribution. Hence, we do not pursue a theoretical comparison between the powers of tests.

### 4·2. *Power analysis*

The following theorem characterizes the conditions under which the power of $T_\alpha$ tends to unity as $n$ grows, under the alternative $H_a^V$ in (14).

THEOREM 3. *Assume that the alignment assumption in* (2) *holds, and that* $\sigma_V^2 = A_1\{1 + o(1)\}$, $\max\{|f(0)|, |g(0)|\} \leq A_2$ *for some positive constants* $A_1$ *and* $A_2$. *Further assume that* $f(\cdot), g(\cdot)$ *have bounded Lipschitz constants. Then for some large scalar* $B_1$ *depending only on* $A_1, A_2$ *and the Lipschitz constants of* $f(\cdot)$ *and* $g(\cdot)$,

$$\inf_{F(X) \in H_a^V(B_1)} \text{pr}(T_\alpha = 1) = 1 - o(1),$$

*where the infimum is taken over all distributions* $F(X)$ *such that* $V\{F(X)\} \in \mathcal{U}(B_1)$.

Similarly, $\widetilde{T}_\alpha$ attains the power tending to unity under the alternative $H_a^U$ in (16).

THEOREM 4. *Suppose that the kernel function* $h(\cdot)$ *in* (3) *is bounded with* $|h(\cdot)| \leq A_3$ *and*

$$m^2 \text{var}_{H_0}\left[E_{H_0}\{h(X_{1,\{1,2\}}, \ldots, X_{m,\{1,2\}}) \mid X_{1,\{1,2\}}\}\right] = \{1 + o(1)\}A_4$$

*for some positive constants* $A_3$ *and* $A_4$. *Then for some large scalar* $B_2$ *depending only on* $A_3, A_4$, *and* $m$,

$$\inf_{F(X) \in H_a^U(B_2)} \text{pr}(\widetilde{T}_\alpha = 1) = 1 - o(1),$$

*where the infimum is taken over all distributions* $F(X)$ *such that* $U\{F(X)\} \in \mathcal{U}(B_2)$.

Here $H_a^V(B_1)$ and $H_a^U(B_2)$ are both sparse alternatives, which can be very close to the null in the sense that all but a small number of entries in $V$ or $U$ can be exactly zero. The above theorems show that the proposed tests are sensitive to small perturbations to the null. Considering the examples discussed in Section 2, Theorems 3 and 4 show that their powers tend to unity under the sparse alternative.

### 4·3. *Optimality*

We now establish the optimality of the proposed tests in the following sense: recall that $T_\alpha$ and $\widetilde{T}_\alpha$ can correctly reject the null provided that at least the magnitude of one entry of $V$ or $U$ is larger than $C(\log d/n)^{1/2}$ for some constant $C$. We show that such a bound is rate-optimal, i.e., the rate of the signal gap, $(\log d/n)^{1/2}$, cannot be further relaxed.

For each $n$, define $\mathcal{T}_\alpha$ to be the set of all measurable size-$\alpha$ tests. In other words, $\mathcal{T}_\alpha \equiv \{T_\alpha : \text{pr}(T_\alpha = 1 \mid H_0) \leq \alpha\}$.

THEOREM 5. *Assume that* $c_0 < 1$ *is a positive constant. Further let* $\beta$ *be a positive constant satisfying that* $\alpha + \beta < 1$. *Under the regime* $\log d/n = o(1)$ *as* $n$ *grows, for large* $n$ *and* $d$, *we have*

$$\inf_{T_\alpha \in \mathcal{T}_\alpha} \sup_{F(X) \in H_a^R(c_0)} \text{pr}(T_\alpha = 0) \geq 1 - \alpha - \beta,$$

*where the supremum is taken over the distribution family* $H_a^R(c_0)$ *defined in* (17).

Theorem 5 shows that any measurable size-$\alpha$ test cannot differentiate the null hypothesis $H_0$ and the sparse alternative when $\max_{j<k}|R_{jk}| \leq c_0(\log d/n)^{1/2}$ for some constant $c_0 < 1$.

The Supplementary Material gives the detailed proof of Theorem 5. It begins with the observation that the family of alternative distributions $H_a^R(c_0)$ includes some Gaussian distributions as a subset, since one can construct a Gaussian distribution given any $R \in \mathcal{U}(c_0)$. Therefore, the supremum over $H_a^R(c_0)$ is no smaller than the supremum over the Gaussian subset. The rest of the proof follows from the general framework in Baraud (2002). In particular, the proof tech-nique is relevant to that used in deriving the lower bound in two-sample covariance tests (Cai et al., 2013).

Due to technical constraints, the alternative considered in Theorem 5 is defined with the Pearson's population correlation matrix $R$. As mentioned in Section 4·1, it is unclear whether there exist equivalences between $V$, $U$, or $R$. Hence, in order to apply Theorem 5 to the alternatives used in Theorems 3 and 4, we make the following assumptions.

*Assumption* 1. When $X$ is Gaussian, the matrices $V$ and $R$ satisfy that, for large $n$ and $d$, $cV_{jk} \leq R_{jk} \leq CV_{jk}$ for $j \neq k \in \{1, \ldots, d\}$, where $c$ and $C$ are two constants.

*Assumption* 2. When $X$ is Gaussian, the matrices $U$ and $R$ satisfy that, for large $n$ and $d$, $cU_{jk} \leq R_{jk} \leq CU_{jk}$ for $j \neq k \in \{1, \ldots, d\}$, where $c$ and $C$ are two constants.

In the proof of Theorem 5, we obtain a lower bound by discussing the Gaussian subset of $H_a^R(c_0)$. This is why we require Assumptions 1 and 2 to hold at least for the Gaussian distributions. Theorem 5 hence justifies the rate-optimality of the proposed tests, summarized as follows.

THEOREM 6. *(a) Suppose that the simple linear rank statistics* $\{V_{jk}, 1 \leq j < k \leq d\}$ *satisfy all conditions in Theorems* 1 *and* 3. *Further suppose that Assumption* 1 *holds. We then have, under the regime* $\log d = o(n^{1/3})$ *as* $n$ *grows, the corresponding size-*$\alpha$ *test* $T_\alpha$ *is rate-optimal. In other words, there exist two constants* $D_1 < D_2$ *such that: (a.i)*

$$\sup_{F(X) \in H_a^V(D_2)} \mathrm{pr}(T_\alpha = 0) = o(1);$$

*(a.ii) for any* $\beta > 0$ *satisfying that* $\alpha + \beta < 1$, *for large* $n$ *and* $d$,

$$\inf_{T_\alpha \in \mathcal{T}_\alpha} \sup_{F(X) \in H_a^V(D_1)} \mathrm{pr}(T_\alpha = 0) \geq 1 - \alpha - \beta.$$

*(b) For all rank-type* $U$-*statistics satisfying the conditions in Theorems* 2 *and* 4, *supposing that Assumption* 2 *holds, then the same rate-optimality property holds.*

As an example, the next corollary justifies the test statistics in Examples 1–4.

COROLLARY 3. *The four test statistics in Examples* 1–4 *are all rate-optimal against the corresponding sparse alternative.*

Corollary 3 is a direct consequence of Theorem 6 and Lemma C8 in the Supplementary Material. Its proof is accordingly omitted.

## 5. NUMERICAL EXPERIMENTS
### 5·1. *Tests*

We compare the performances of our proposed tests with competitors on various synthetic data sets in both low-dimensional and high-dimensional settings. We provide additional numerical results in the Supplementary Material, which includes comparisons to other tests of $H_0$, testings with simulation-based rejection thresholds, and an application.

We propose a test based on the Spearman's rho statistic, outlined in Example 1, using the result in Theorem 1. We propose another test based on the Kendall's tau statistic, outlined in Example 2, using the result in Theorem 2. We use the theoretical rejection threshold $q_\alpha$ in (9) for both tests. We will refer to these tests as the Spearman test and the Kendall test below.

We consider the test of Zhou (2007), which rejects the null if

$$n \max_{j<k} r_{jk}^2 - 4\log d + \log\log d \geq q_\alpha, \tag{18}$$

where $q_\alpha$ is the theoretical threshold defined in (9).

As another competitor, we consider the test of Mao (2014), which rejects the null hypothesis if

$$\left\{\sum_{j<k} \frac{r_{jk}^2}{1-r_{jk}^2} - \frac{d(d-1)}{2(n-4)}\right\} \left\{\frac{(n-4)^2(n-6)}{d(d-1)(n-3)}\right\}^{1/2} \geq \Phi^{-1}(1-\alpha), \tag{19}$$

where $\Phi^{-1}(\cdot)$ is the quantile function of the standard Gaussian. The test statistic in Mao (2014) has guaranteed size control only under normality.

The last two competitors are the kernel-based tests of Reddi & Póczos (2013) and Póczos et al. (2012). Briefly speaking, Reddi & Póczos (2013) propose to calculate the Hilbert–Schmidt norm of the normalized cross-covariance operators after a copula transformation, and Póczos et al. (2012) propose to use the estimated maximum mean discrepancy after a copula transformation. In both kernel-based tests, we use the Gaussian kernel with standard deviation as the median distance heuristic as in Reddi & Póczos (2013) and Póczos et al. (2012). We use simulation to determine the rejection thresholds for both tests, since the null distribution of $F(X)$ is uniform. Although a theoretical rejection threshold is proposed in Póczos et al. (2012), it becomes too conservative in high dimensions.

In summary, we apply six tests in the numerical experiments, namely, the Spearman test outlined in Example 1, the Kendall test outlined in Example 2, the test of Zhou (2007), the test of Mao (2014), the test of Reddi & Póczos (2013), and the test of Póczos et al. (2012). In the following experiments, we set the nominal significance level to be $\alpha = 0 \cdot 05$ for all tests.

We have further compared the performance of our proposed procedures with two more rank-based statistics (Mao, 2016; Leung & Drton, 2017). Due to space limit, these additional results are put in the Supplementary Material.

### 5·2. *Synthetic data analysis*

We now provide size and power comparisons among the competing tests introduced in Section 5·1. In this simulation, we generate synthetic data from five different types of distributions: the Gaussian distribution, the light-tailed Gaussian copula, the heavy-tailed Gaussian copula, the multivariate $t$ distribution, and the multivariate exponential distribution. To evaluate the sizes of the tests, we generate data from the five types of distributions under the null, where all entries in $X$ are mutually independent. For evaluating the powers of the tests, we generate different sets of data from the five types of distributions under sparse alternatives. For instance, for the Gaussian distribution, we draw our data from $N_d(0, I_d)$ to evaluate the size, and generate data from $N_d(0, R^*)$ to evaluate the power. Here $R^* \in \mathcal{R}^{d \times d}$ is a positive definitive matrix, whose off-diagonal entries are all zero except for eight randomly chosen entries. We defer details of the data generating mechanism to the Supplementary Material.

In summary, we generate data from ten distributions: one under the null and one under the sparse alternative for each of the five types. For each distribution, we draw $n$ independent replicates of the $d$-dimensional random vector $X \in \mathcal{R}^d$. To examine the effects of increasing sample sizes and dimensions, let the sample size $n$ be 60 or 100, and the dimension $d$ be 50, 200, or 800. Results from 5,000 simulated data sets are given in Figure 1.

Under the Gaussian distribution, most tests can effectively control the size under all combinations of $n$ and $d$. Zhou (2007) and our proposals attain higher power than the others. In contrast, the test of Mao (2014) has relatively lower power. This is as expected because, under the sparse alternative, the correlation matrix has only eight non-zero entries. By averaging over all entries in the correlation matrix, the test of Mao (2014) is less sensitive to the sparse alternative. Due to similar reasons, the test of Póczos et al. (2012) has low power against the sparse alternative. The test of Reddi & Póczos (2013) has decreasing power as $d$ grows, which corroborates the findings in Ramdas et al. (2015).

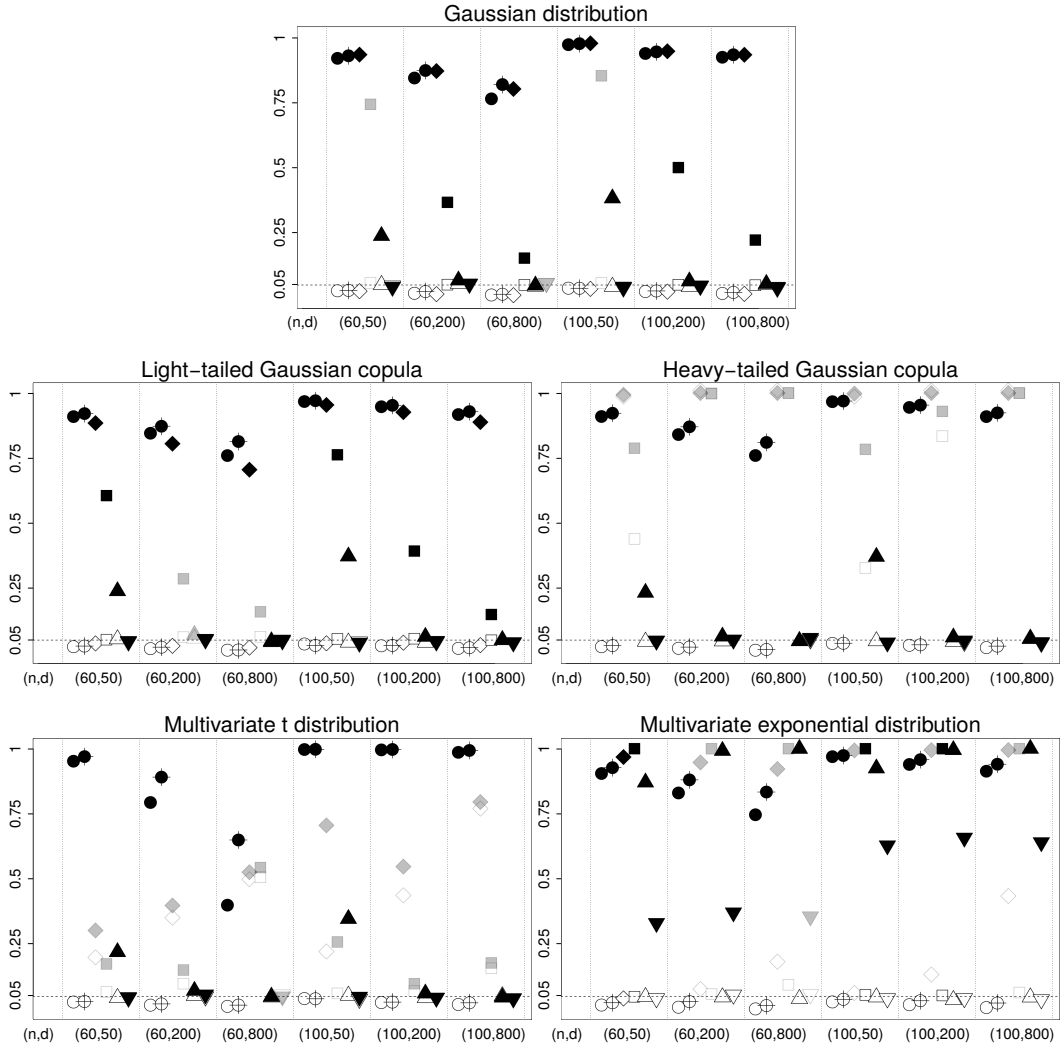

Fig. 1: Empirical sizes and powers of six tests under five types of distributions averaged over 5,000 replicates at the $0 \cdot 05$ nominal significance level, shown by the horizontal dashed line. The y-axis represents the proportion of rejected tests in the 5,000 replicates. The vertical dotted lines separate six different data generating schemes, where the pair of sample size and dimension, i.e., $(n, d)$, range from $(60, 50)$ to $(100, 800)$. Only $(60, 50)$ and $(100, 50)$ are in the low-dimensional setting. The six tests considered in the simulation are the Spearman test (○ and ●), the Kendall test (⊕ and ⬥), the test of Zhou (2007) (◇ and ◆), the test of Mao (2014) (□ and ■), the test of Reddi & Póczos (2013) (△ and ▲), and the test of Póczos et al. (2012) (▽ and ▼). We use hollow shapes to represent the empirical sizes under the null (e.g., ⊕), and solid shapes to represent the empirical powers under the alternative (e.g., ⬥). A gray symbol indicates that the corresponding test fails to control the size at the $0 \cdot 05$ nominal significance level in the corresponding null model (e.g., △ and ▲). In this simulation, we say that a test fails to control the size at $0 \cdot 05$ if its empirical size is larger than $0 \cdot 05 + 1 \cdot 96(0 \cdot 05 \times 0 \cdot 95/5000)^{1/2}$.

Under the light-tailed Gaussian copula, the performances of all six tests are similar to their performances under the Gaussian distribution. Notably, our proposed tests achieve higher power than the test of Zhou (2007), especially when the ratio $d/n$ is large.

Under the heavy-tailed Gaussian copula, only our proposed tests correctly control the size while attaining high power. Zhou (2007) and Mao (2014) have very high power under the alternative, but their sizes are severely inflated under the null. The kernel-based tests (Reddi & Póczos, 2013; Póczos et al., 2012) correctly control the size but have very low power.

Under the multivariate $t$ distribution, the performances are similar to those under the heavy-tailed Gaussian copula.

Under the multivariate exponential distribution, our proposals and the test of Reddi & Póczos (2013) achieve high power while correctly controlling the size across all settings. The tests of Mao (2014) and Zhou (2007) fail to control the size in high dimensions. The test of Póczos et al. (2012) has low power compared to the others.

In summary, the proposed tests correctly control the size and achieve high power across all types of distributions regardless of the sample size and dimension. The Kendall test performs slightly better than the Spearman test in terms of power. This is consistent with the observations in Woodworth (1970), which shows that Kendall's tau is asymptotically more powerful than Spearman's rho in testing independence in terms of having Bahadur efficiency bounded in $(1, 1 \cdot 05]$ under the Gaussian distribution. For relationship between Spearman's rho and Kendall's tau under different alternatives, we refer to Fredricks & Nelsen (2007b) for details. The performances of Zhou (2007) and Mao (2014) are severely influenced by the data structure, and cannot effectively control the size under heavy-tailed distributions. This is as expected because the validity of Mao (2014)'s test relies heavily on the Gaussian assumption, and the performance of Zhou (2007) is related to the moments. The kernel-based tests (Reddi & Póczos, 2013; Póczos et al., 2012) control the size correctly in most cases, but their power suffers in high dimensions, as observed in Ramdas et al. (2015).

## 6. DISCUSSION

The regression constants $\{c_{ni}\}_{i=1}^{n}$, the score function $g(\cdot)$ in (1), and the kernel function $h(\cdot)$ in (3) are assumed to be identical across different pairs of entries, which can be straightforwardly relaxed. We do not pursue this direction merely for clarity of presentation.

The problem studied in this manuscript is related to one-sample and two-sample tests of equality of covariance or correlation matrices and sphericity tests in high dimensions. There exist extensive studies along this line of research. See, among others, Ledoit & Wolf (2002), Chen et al. (2010), Srivastava & Yanagihara (2010), Fisher et al. (2010), Li & Chen (2012), Fisher (2012), Zhang et al. (2013), Cai et al. (2013), Han et al. (2017). For equity and sphericity tests, existing methods mostly focus on the Pearson's sample covariance matrix. Corresponding tests are based on statistics characterizing the difference between two-sample covariance matrices under different norms, such as the Frobenious norm or the maximum norm. As an alternative, Zou et al. (2014) propose a sphericity test using the multivariate signs. However, the theoretical results in their paper are valid only under the regime $d = O(n^2)$.

Testing equality of covariance or correlation matrices is challenging since the random variables are not mutually independent. In the Supplementary Material, focusing on the one-sample test, we test the bandedness of the latent correlation matrix under the semi-parametric Gaussian copula model. We show that the test built on Kendall's tau statistic can asymptotically control the size, and is rate-optimal against the sparse alternative. More details are relegated to the Supplementary Material.

ACKNOWLEDGEMENT

We thank for helpful discussions with Cheng Zhou and Brian Caffo. We thank the editor, associate editor, and two referees for their valuable suggestions. This research was supported in part by the U.S. National Science Foundation and National Institutes of Health.

SUPPLEMENTARY MATERIAL

Supplementary material available at *Biometrika* online includes discussion of simulation-based rejection thresholds, generalizations of the proposed tests to other structural testing problems, additional numerical results, and proofs for the theoretical claims in the main paper.

# Supplementary Material for Distribution-Free Tests of Independence in High Dimensions

By FANG HAN

*Department of Statistics, University of Washington, Box 354322, Seattle, Washington 98195, U.S.A.*

fanghan@uw.edu

SHIZHE CHEN

*Department of Statistics, Columbia University, Room 1005 SSW, MC 4690, 1255 Amsterdam Avenue, New York, New York 10027, U.S.A.*

shizhe.chen@gmail.com

AND HAN LIU

*Department of Operations Research and Financial Engineering, Princeton University, Sherrerd Hall, Charlton Street, Princeton, New Jersy 08544, U.S.A.*

hanliu@princeton.edu

## A. ADDITIONAL RESULTS

### A·1. *Overview*

The theory in the main paper employs techniques that can be easily generalized to other problems such as structural testings. In this section, we discuss three additional results that are of interest. In particular, Section A·2 studies the approximation of the exact distributions of the test statistics proposed in the main paper, and we consider the problems of testing $m$-dependence and homogeneity in Sections A·3 and A·4.

### A·2. *Approximation to the exact distributions*

Theorems 1 and 2 in the main paper show that the proposed test statistics $L_n$ and $\widetilde{L}_n$ converge weakly to a Gumbel distribution. The next theorem characterizes the convergence rates for $L_n$ and $\widetilde{L}_n$.

THEOREM A1. *For all rank-type $U$-statistics, under the conditions in Theorem* 2 *and that* $\log d = o(n^{1/3})$, *we have*

$$\Big|\mathrm{pr}\Big(\frac{n\widetilde{L}_n^2}{\sigma_U^2} - 4\log d + \log\log d \leq y\Big) - \exp\Big\{-(8\pi)^{-1/2}\exp\Big(-\frac{y}{2}\Big)\Big\}\Big| = O_y\Big\{\frac{(\log d)^{3/2}}{n^{1/2}} + \frac{1}{(\log d)^{3/2}}\Big\}.$$

*For all simple linear rank statistics, if conditions in Theorem* 1 *hold and* $\log d = O(n^{1/3-\epsilon})$ *for some constant* $\epsilon \in (0, 1/3)$, *we have*

$$\Big|\mathrm{pr}\Big(\frac{nL_n^2}{\sigma_V^2} - 4\log d + \log\log d \leq y\Big) - \exp\Big\{-(8\pi)^{-1/2}\exp\Big(-\frac{y}{2}\Big)\Big\}\Big|$$
$$= O_y\Big\{\frac{(\log d)^{3/2}}{n^{1/2}} + \frac{1}{(\log d)^{3/2}} + \frac{(\log d)^{1/2}}{n^{1/6}}\Big\}.$$

Theorem A1 shows two points. (i) When $\log d \asymp n^\kappa$ for some $\kappa < 1/3$, the proposed tests based on simple linear rank statistics and rank-type $U$-statistics achieve polynomial rates of convergence. Compared to tests based on the rank-type $U$-statistics, the tests based on simple linear rank statistics lose an extra $O\{(\log d)^{1/2} n^{-1/6}\}$ term in the rate of convergence, due to approximating the population ranks using the empirical ranks. Check the proof of moderate deviation in Lemma C5 for more details. (ii) When $d \asymp n^C$ for some $C \in (0, \infty)$, Theorem A1 only guarantees an $O\{(\log n)^{-3/2}\}$ rate of convergence.

We will show that the convergence rate can be accelerated by approximating the exact distributions of the test statistics. Under $H_0$ in the main paper, $\{V_{jk}, j < k\}$ and $\{U_{jk}, j < k\}$ are independent and only depend on the relative ranks $\{R_{ni}^{jk}, i = 1, \ldots, n, j < k\}$, which are uniformly distributed under permutations on $\{1, \ldots, n\}$. Therefore, we can conduct simulations to approximate the exact distributions of $\{V_{jk}, j < k\}$ and $\{U_{jk}, j < k\}$, respectively.

Specifically, for $i = 1, \ldots, M$, we generate $X_{\cdot,\cdot}^{(i)} \in \mathcal{R}^{n \times d}$ as an $n \times d$ matrix with all entries independently drawn from a standard normal distribution, which yield simple linear rank statistics $\{V_{jk}^{(i)}, j < k\}$ and the rank-type $U$-statistics $\{U_{jk}^{(i)}, j < k\}$. Next, we calculate the values of $n(L_n^{(i)})^2/\sigma_V^2 - 4\log d + \log \log d$ and $n(\widetilde{L}_n^{(i)})^2/\sigma_U^2 - 4\log d + \log \log d$. Here $L_n^{(i)}$ and $\widetilde{L}_n^{(i)}$ are the extreme-value statistics based on $\{V_{jk}^{(i)}, j < k\}$ and $\{U_{jk}^{(i)}, j < k\}$, respectively. Let $\widehat{F}_{n,d;M}^V(\cdot)$ and $\widehat{F}_{n,d;M}^U(\cdot)$ be the empirical distributions, and let $F_{n,d}^V(\cdot)$ and $F_{n,d}^U(\cdot)$ be their population counterparts.

The Dvoretzky-Kiefer-Wolfowitz inequality (Dvoretzky et al., 1956; Massart, 1990) guarantees, for each pair of $(n, d)$,

$$\begin{aligned}\text{pr}\left\{\sup_{x \in \mathcal{R}} |\widehat{F}_{n,d;M}^V(x) - F_{n,d}^V(x)| > \left(\frac{\log M}{M}\right)^{1/2}\right\} &\leq \frac{2}{M^2}, \\ \text{pr}\left\{\sup_{x \in \mathcal{R}} |\widehat{F}_{n,d;M}^U(x) - F_{n,d}^U(x)| > \left(\frac{\log M}{M}\right)^{1/2}\right\} &\leq \frac{2}{M^2}.\end{aligned} \quad \text{(A1)}$$

We replace $q_\alpha$ in (8) using $\widehat{q}_{\alpha;n,d}^V$ and $\widehat{q}_{\alpha;n,d}^U$, which are the $1 - \alpha$ quantiles of $\widehat{F}_{n,d;M}^V(\cdot)$ and $\widehat{F}_{n,d;M}^U(\cdot)$

$$\widehat{q}_{\alpha;n,d}^V \equiv \inf\{x : \widehat{F}_{n,d;M}^V(x) \geq 1 - \alpha\}, \quad \widehat{q}_{\alpha;n,d}^U \equiv \inf\{x : \widehat{F}_{n,d;M}^U(x) \geq 1 - \alpha\}.$$

We refer to the tests using the simulation-based thresholds $\widehat{q}_{\alpha;n,d}^V$ and $\widehat{q}_{\alpha;n,d}^U$ as the exact tests.

Using (A1), we have the next theorem that guarantees the asymptotic control of sizes.

THEOREM A2. *Under $H_0$, simple linear rank statistics satisfy that, for each pair of $(n, d)$, with probability no smaller than $1 - 2/M^2$, we have*

$$\sup_{\alpha \in [0,1]} \left|\text{pr}\left(\frac{nL_n^2}{\sigma_V^2} - 4\log d + \log\log d \geq \widehat{q}_{\alpha;n,d}^V \mid \{X_{\cdot,\cdot}^{(i)}\}_{i=1}^M\right) - \{1 - \widehat{F}_{n,d;M}^U(\widehat{q}_{\alpha;n,d}^V)\}\right| \leq \left(\frac{\log M}{M}\right)^{1/2}.$$

*The same inequality also applies to the rank-type $U$-statistics. Moreover, as $n$ and $d$ grow, $\widehat{q}_{\alpha;n,d}^V$ and $\widehat{q}_{\alpha;n,d}^U$ are both consistent estimators of $q_\alpha$ in (9) as $M = M_n$ grows with $n$.*

Theorem A2 shows that, with high probability, we can have arbitrarily fast convergence rates to the above intermediate approximation by setting the number of simulations $M$ large enough. Typically, it is much faster than the rate $O\{(\log n)^{5/2}/n^{1/2}\}$ derived in Liu et al. (2008). On the other hand, to attain this arbitrarily fast rate of convergence, we need to conduct $M$ simulations for estimating the threshold value. This increases the computational burden compared to the tests in (8). For the test of $m$-dependence, which we shall introduce in Section A·3, it is impossible to simulate the null exact distribution and we stick to the test in (A2).

### A·3. *Test of $m$-dependence*

A random vector $X = (X_1, \ldots, X_d)^{\text{T}} \in \mathcal{R}^d$ follows a Gaussian copula distribution if and only if $\{F_1(X_1), F_2(X_2), \ldots, F_d(X_d)\}^{\text{T}}$ distributes the same as $\{\Phi(Z_1), \ldots, \Phi(Z_d)\}^{\text{T}}$, where $F_1, \ldots, F_d$ are the marginal distribution functions of $X_1, \ldots, X_d$, $\Phi(\cdot)$ represents the distribution function of the standard Gaussian, and $Z = (Z_1, \ldots, Z_d)^{\text{T}} \sim N_d(0, \Sigma^0)$ with diagonal entries of $\Sigma^0$ equal 1. The Gaussian copula family includes the Gaussian, and is a semi-parametric one since the marginal distributions of $X$ are unspecified. We refer to $\Sigma^0$ as the latent correlation matrix of $X$. As in the main paper, we only consider continuous $X$ for avoiding possible ties.

We aim at testing the null hypothesis $A_0 : \Sigma^0_{jk} = 0$, for all $|j - k| \geq m$. Because $X$ is assumed to be a Gaussian copula, the dependence structure among $\{X_1, \ldots, X_d\}$ is fully encoded in $\Sigma^0$. Therefore, testing $A_0$ is equivalent to testing $m$-dependence among entries of $X$, i.e., $X_j$ is independent of $X_k$, for all $|j - k| \geq m$.

Cai & Jiang (2011) first consider the problem of testing $A_0$ in high dimensions on Gaussian data. Later, the result is extended to non-Gaussian data under a moment assumption (Shao & Zhou, 2014). In this section, we show that the moment assumption can be utterly relaxed by resorting to the rank-based statistics.

For testing $A_0$, instead of resorting to the Pearson's correlation coefficients as in Cai & Jiang (2011) and Shao & Zhou (2014), we use Kendall's tau correlation coefficients $\{\tau_{jk}, 1 \leq j < k \leq d\}$ introduced in Example 2 in the main paper. It is well known that Kendall's tau is irrelevant to the marginal distributions of $X$ (Nelsen, 1999). Accordingly, within the Gaussian copula family, Kendall's tau is a more natural measure of dependence than Pearson's correlation coefficient. Moreover, it is known from Lemma C8 that, under the Gaussian copula family, we have $\Sigma^0_{jk} = \sin(\tau^0_{jk} \pi/2)$, where $\tau^0_{jk} \equiv E(\tau_{jk})$. Therefore, within the Gaussian copula family, testing $A_0$ is equivalent to testing $\tau^0_{jk} = 0$ for all $|j - k| \geq m$. We hence propose the following test statistic

$$T^\tau_{\alpha, m} \equiv I\left\{\frac{9n}{4}(L^\tau_{n,m})^2 - 4\log d + \log\log d \geq q_\alpha\right\}, \tag{A2}$$

where $q_\alpha$ is introduced in (9) in the main paper and the extreme-value statistic $L^\tau_{n,m} \equiv \max_{|j-k| \geq m} |\tau_{jk}|$. $L^\tau_{n,m}$ is an extreme-value statistic similar to $L^\tau_n$ in the main paper. We expect $L^\tau_{n,m}$ to have similar null limiting distribution as $L^\tau_n$ given proper conditions on $m$. We reject $A_0$ if and only if $T^\tau_{\alpha, m} = 1$.

The following theorem justifies the test $T^\tau_{\alpha, m}$ for a fixed nominal significance level $\alpha$.

THEOREM A3. *Suppose that $\log d = o(n^{1/3})$ as $n$ grows, $m = o(d^c)$ for any $c > 0$, and for some constant $\delta \in (0, 1)$,*

$$\text{card}\bigl[\{1 \leq j \leq d : |\Sigma^0_{jk}| > 1 - \delta \text{ for some } 1 \leq k \leq d \text{ and } j \neq k\}\bigr] = o(d).$$

*Provided that $X$ is continuous and distributes as a Gaussian copula, under $A_0$, we have, for any $y \in \mathcal{R}$,*

$$\left|\text{pr}\left\{\frac{9n}{4}(L^\tau_{n,m})^2 - 4\log d + \log\log d \leq y\right\} - \exp\left\{-(8\pi)^{-1/2}\exp\left(-\frac{y}{2}\right)\right\}\right| = o_y(1).$$

*Accordingly, the test $T^\tau_{\alpha, m}$ can asymptotically control the size as $n$ and $d$ grow, i.e.,*

$$\text{pr}(T^\tau_{\alpha, m} = 1 \mid A_0) = \alpha + o(1).$$

*Remark* A1. The proof of the theorem shows that the assumption, $m = o(d^c)$ for any $c > 0$, can be easily relaxed. Specifically, we only require $m = o(d^{\epsilon(\delta)})$ for a small enough constant $\epsilon(\delta)$ depending on $\delta$. This can be verified by checking Equation (C19), and Equation (68) in Cai & Jiang (2011).

Similar to the power analysis in Section 4·2 in the main paper, we study the power of the test $T^\tau_{\alpha, m}$ against a sparse alternative. To this end, consider the following set of matrices

$$\mathcal{U}_m(c) \equiv \Bigl\{M \in \mathcal{R}^{d \times d} : \text{diag}(M) = I_d, M = M^{\mathrm{T}}, \max_{|j-k| \geq m} |M_{jk}| \geq c(\log d/n)^{1/2}\Bigr\}.$$

The following theorem shows, for the Gaussian copula family, as long as the latent correlation matrix $\Sigma^0 \in \mathcal{U}_m(C)$ for some large constant $C$, the power of the proposed test tends to one.

THEOREM A4. *Suppose that we observe $n$ independent observations of a $d$-dimensional random vector $X = (X_1, \ldots, X_d)^{\mathrm{T}}$ following a Gaussian copula with the latent correlation matrix $\Sigma^0$. Then, there* *exists some large constant $D_3$ such that*

$$\sup_{\Sigma^0 \in \mathcal{U}_m(D_3)} \mathrm{pr}(T^\tau_{\alpha,m} = 1) = 1 - o(1),$$

*as $n$ and $d$ grow. Here the supremum is taken over the Gaussian copula family such that $\Sigma^0 \in \mathcal{U}_m(D_3)$.*

We derive Theorem A4 using a similar technique as in the proof of Theorem 3. The proof is thus omitted.

We then turn to study the optimality of $T^\tau_{\alpha,m}$. In testing $A_0$, for each $n$, we define $\mathcal{T}_{\alpha,m}$ to be the set of all measurable size-$\alpha$ tests $T_{\alpha,m}$ such that $\mathrm{pr}(T_{\alpha,m} = 1 \mid A_0) \leq \alpha$. The following theorem gives the detection lower bound in differentiating the null hypothesis and the sparse alternative.

THEOREM A5. *Assume that there exists a positive constant $c'_0 < 1$, $\log d = o(n)$ as $n$ grows, and $m = o(d^c)$ for any $c > 0$. Let $\beta$ be a positive constant satisfying that $\alpha + \beta < 1$. For all large enough $n$ and $d$, we have*

$$\inf_{T_{\alpha,m} \in \mathcal{T}_{\alpha,m}} \sup_{\Sigma^0 \in \mathcal{U}_m(c'_0)} \mathrm{pr}(T_{\alpha,m} = 0) \geq 1 - \alpha - \beta,$$

*where the supremum is taken over any distribution family such that $\Sigma^0 \in \mathcal{U}_m(c'_0)$.*

Therefore, we conclude that $T^\tau_{\alpha,m}$ is rate-optimal in testing the null hypothesis $A_0$ against the sparse alternative in the main paper.

For any constant $c > 0$, the matrix set $\mathcal{U}(c)$ defined in (13) in the main paper includes $\mathcal{U}_m(c)$. Accordingly, the lower bound derived in Section 4·3 cannot be trivially exploited to derive the lower bound for testing the bandedness of $\Sigma^0$. However, using the fact that $m = o(d^c)$ for any $c > 0$, we can find the lower bound for testing $A_0$ via designing a similar set of parameters as in the proof of Theorem 5.

### A·4. *Test of homogeneity*

Let $X_{1,\cdot}, \ldots, X_{n,\cdot} \in \mathcal{R}^d$ be $n$ independent but not necessarily identically distributed random vectors with $X_{i,\cdot} = (X_{i,1}, \ldots, X_{i,d})^\mathrm{T}$ for $i = 1, \ldots, n$. We aim at testing $B_0 : X_{1,\cdot}, \ldots, X_{n,\cdot}$ are identically distributed. Testing $B_0$ is of fundamental interest in many statistical fields.

It is generally very complicated to test homogeneity in high dimensions. The works in this field are very limited and most of the existed ones reduce it to equity tests of two-sample means and covariance matrices. Bai & Saranadasa (1996), Srivastava & Du (2008), Chen & Qin (2010), and Cai et al. (2014) consider comparing the means of two high-dimensional Gaussian vectors with unknown covariance matrices, and Chen et al. (2010) and Cai et al. (2014) develop tests of equity of two covariance matrices.

We consider a simplified version of $B_0$: the entries in each $X_{i,\cdot}$ are mutually independent. In this simplified setting, we reduce the test of $B_0$ to the test that $X_{1,j}, X_{2,j}, \ldots, X_{n,j}$ are identically distributed for any $j \in \{1, \ldots, d\}$. For each $j$, we test the homogeneity using a rank-based test statistic. We then formulate an extreme-value statistic by combining the $d$ separate rank-based test statistics.

In details, let $H_n$ be an extreme-value statistic summarizing the $d$ separate rank-based test statistics: $H_n \equiv \max_{j \in \{1,\ldots,d\}} |h_j|$, where

$$h_j \equiv \frac{2}{n(n-1)} \sum_{i<i'} \mathrm{sign}(X_{i',j} - X_{i,j}) \quad (j = 1, \ldots, d).$$

Here $h_j$ is an rank-based statistic counting the number of inequalities $X_{i',j} > X_{i,j}$ across all pairs $i < i'$. Mann (1945) is the first to introduce the test statistic $h_j$ for testing homogeneity. Mann (1945) characterizes the sufficient conditions for $h_j$ to be consistent and unbiased, and shows that this statistic is powerful against a trend alternative that will be introduced later. We refer to Kendall & Stuart (1961) for more discussion on the rationale of using $h_j$ for testing homogeneity. For testing $B_0$, we propose the following statistic based on $H_n$:

$$T^h_\alpha \equiv I\Big(\frac{9n}{4} H_n^2 - 2\log d + \log\log d \geq \widetilde{q}_\alpha\Big),$$

where $\widetilde{q}_\alpha \equiv -\log \pi - 2\log\log(1-\alpha)^{-1}$ is the $1-\alpha$ quantile of the Gumbel distribution with the distribution function $\exp\{-\pi^{-1/2}\exp(-y/2)\}$.

Next, we justify that the test $T_\alpha^h$ controls the size properly. Under $B_0$, we have $X_{1,j},\ldots,X_{n,j}$ are identically distributed and hence the distribution of $\mathrm{sign}(X_{i',j} - X_{i,j})$ should be centered around zero, and the ranks of $X_{1,j},\ldots,X_{n,j}$ are uniformly sampled from the set of all permutations of $\{1,\ldots,n\}$. Accordingly, $h_j$ is identically distributed to Kendall's tau statistic under $H_0$ in the main paper. Therefore, using Example 2, we derive $E_{B_0}(h_j) = 0$ and

$$\mathrm{var}_{B_0}(h_j) = \frac{2(2n+5)}{9n(n-1)} = \frac{4}{9n}\{1+o(1)\},$$

and the limiting distribution of $H_n$ shall resemble that of Kendall's tau. Specifically, the following theorem provides the null limiting distribution of $H_n$.

THEOREM A6. *Suppose that $\log d = o(n^{1/3})$ as $n$ grows. Under $B_0$, we have, for any $y \in \mathcal{R}$,*

$$\left|\mathrm{pr}\left(\frac{9n}{4}H_n^2 - 2\log d + \log\log d\right) - \exp\left\{-\pi^{-1/2}\exp\left(-\frac{y}{2}\right)\right\}\right| = o_y(1).$$

*Accordingly, the test $T_\alpha^h$ can asymptotically control the size as $n$ and $d$ grow, i.e.,*

$$\mathrm{pr}(T_\alpha^h = 1 \mid B_0) = \alpha + o(1).$$

It is worth noting that, similar to Corollary 1 in the main paper, Theorem A6 holds without any distributional assumption on $X_{1,\cdot},\ldots,X_{n,\cdot}$.

We then study the power of the proposed test. We consider a particular trend alternative; that is, for at least one entry $j \in \{1,\ldots,d\}$, the mean of $X_{i,j}$ is a linear function of $i$ for a certain entry $j \in \{1,\ldots,d\}$, i.e., $B_1$: there exists some $j \in \{1,\ldots,d\}$ such that $E(X_{i,j}) = \beta_0 + \beta_1 i/n$ with $\mathrm{var}(X_{i,j}) = \sigma^2$, for $i = 1,\ldots,n$ and $\beta_0, \beta_1, \sigma^2 \in \mathcal{R}$. Under $B_1$, the variance $\sigma^2$ is identical across samples while the means are monotonically increasing or decreasing with respect to the label $i$. Such an alternative is of interest in areas including quality control, finance, and longitudinal data analysis. For instance, in quality control we are interested in inspecting whether machines keep performing well. One alternative of interest is: at least one machine's performance keeps descending.

Under $B_1$, consider the following set of real numbers $(a_1, a_2)$:

$$\mathcal{B}(c) \equiv \left\{(a_1, a_2) : |a_1|/a_2 \geq c(\log d/n)^{1/2}, a_2 > 0\right\}.$$

The following theorem shows that, uniformly over the alternative hypothesis set $\mathcal{B}(C)$, for some large enough constant $C > 0$, the power of the proposed test tends to unity as $n$ grows.

THEOREM A7. *Suppose that there exists at least one entry $j \in \{1,\ldots,d\}$ satisfying $B_1$ with parameters of interest $(\beta_1, \sigma)$. Moreover, for $i = 1,\ldots,n$, the density function $p_{ij}(\cdot)$ of $\{X_{i,j} - E(X_{i,j})\}/\{\mathrm{var}(X_{i,j})\}^{1/2}$ is identical to some density function $p(\cdot)$, which satisfies that*

$$p(x) \geq D_4 > 0 \quad \text{for all } x \in [-M, M], \tag{A3}$$

*for some constant $M > 0$. Then there exists some large scalar $D_5$ only depending on $D_4$ and $M$ such that*

$$\sup_{(\beta_1, \sigma) \in \mathcal{B}(D_5)} \mathrm{pr}(T_\alpha^h = 0) = o(1).$$

In the following we show that the detection boundary $|\beta_1|/\sigma \geq C(\log d/n)^{1/2}$ is rate-optimal. We define $\mathcal{T}_\alpha^h$ to be the set of all measurable size-$\alpha$ tests $T_\alpha^h$ satisfying

$$\mathrm{pr}(T_\alpha^h = 1 \mid B_0) \leq \alpha.$$

The following theorem shows that the proposed test is rate-optimal against the trend alternative $B_1$.

THEOREM A8. *Assume that there exists a constant $c_0'' < 3^{1/2}$, $\log d/n = o(1)$ as $n$ grows. Let $\beta$ be a positive constant satisfying that $\alpha + \beta < 1$. For all large enough $n, d$, we have*

$$\inf_{T_\alpha^h \in \mathcal{T}_\alpha^h} \sup_{(\beta_1, \sigma) \in \mathcal{B}(c_0'')} \mathrm{pr}(T_\alpha^h = 0) \geq 1 - \alpha - \beta,$$

*where $\mathcal{T}_\alpha^h$ represents the family of measurable size-$\alpha$ tests under $B_0$, and the supremum is taken over any distribution family of $X_{1,\cdot}, \ldots, X_{n,\cdot}$ satisfying $B_1$.*

It is straightforward that, when $X_{1,\cdot}, \ldots, X_{n,\cdot}$ are normally distributed, Equation (A3) in Theorem A7 is satisfied. Accordingly, combining Theorems A6, A7, and A8 concludes that $T_\alpha^h$ is rate-optimal in testing the null hypothesis $B_0$ against the trend alternative $B_1$.

## B. Additional numerical experiments

### B·1. *Overview*

In this section, we conduct additional numerical experiments to further explore the properties of our proposals. In Section B·2, we provide details of the data generating mechanism in Section 5·2 in the main paper. In Section B·3, we compare our tests to recent proposals by Mao (2016) and Leung & Drton (2017). In Section B·4, we investigate the performance of the approximation proposal in Section A·2. And finally, we apply our proposals on a real data set in Section B·5.

### B·2. *Data generating mechanism*

We now explain in detail the null distributions and alternative distributions used in Section 5·2 in the main paper.

For the Gaussian distribution, we generate data from $X \sim N_d(0, I_d)$ under the null, and $X \sim N_d(0, R^*)$ under the sparse alternative. Here $R^*$ is generated as follows: consider a random matrix $\Delta \in \mathcal{R}^{d \times d}$ with eight nonzero entries. We select the locations of four nonzero entries randomly from the upper triangle of $\Delta$, each with a magnitude randomly drawn from the uniform distribution in $[0, 1]$. The other four nonzero entries in the lower triangle are determined by symmetry. Finally, to ensure positivity, $R^* \equiv I_d + \Delta + \delta I_d$, where $\delta = \{-\lambda_{\min}(I_d + \Delta) + 0 \cdot 05\} I\{\lambda_{\min}(I_d + \Delta) \leq 0\}$.

For the light-tailed Gaussian copula, we draw data as $X_j = Z_j^{1/3}$ for $j = 1, \ldots, d$ in both the null and alternative distributions. Under the null, $Z = (Z_1, \ldots, Z_d)^{\mathrm{T}} \sim N_d(0, I_d)$, and under the alternative, $Z = (Z_1, \ldots, Z_d)^{\mathrm{T}} \sim N_d(0, R^*)$.

For the heavy-tailed Gaussian copula, we draw data as $X_j = Z_j^3$ for $j = 1, \ldots, d$. Under the null, $Z = (Z_1, \ldots, Z_d)^{\mathrm{T}} \sim N_d(0, I_d)$, and under the alternative, $Z = (Z_1, \ldots, Z_d)^{\mathrm{T}} \sim N_d(0, R^*)$.

For the multivariate $t$ distribution, we generate $X_1, \ldots, X_d$ independently from a univariate $t$ distribution with degree of freedom three under the null distribution, and we generate data from a multivariate $t$ distribution with the covariance matrix $R^*$ and degree of freedom three under the alternative distribution.

For the multivariate exponential distribution, we draw $X_j, j = 1, \ldots, d$ from independent exponential distributions of rate $0 \cdot 25$ under the null distribution, and from a multivariate distribution, where, for each $j = 1, \ldots, d$, $X_j$ conditioned on $X_{-j}$ follows an exponential distribution of rate $0 \cdot 25 + R^*_{j,-j} X_{-j}$. Here $R^*_{j,-j}$ denotes the $j$th row of $R$ without the diagonal element, and $X_{-j}$ denotes the vector $X$ without the $j$th entry.

### B·3. *Additional comparisons*

Mao (2016) and Leung & Drton (2017) study the problem of testing $H_0$ using statistics based on the sums of rank correlations. Mao (2016) proposes a test based on Spearman's rho statistics

$$S = \sigma_{nd}^{-1} \left\{ \sum_{j=2}^{d} \sum_{k=1}^{j-1} \rho_{jk}^2 - \frac{d(d-1)}{2(n-1)} \right\}, \tag{B1}$$

where $\sigma_{nd}^2 \equiv \{d(d-1)(25n^3 - 57n^2 - 40n + 108)\}/\{25(n-1)^3 n(n+1)\}$. Mao (2016) shows that $S$ converges in distribution to the standard normal as $n$ and $d$ grow. Leung & Drton (2017) study a similar statistics

$$T = \frac{n}{d} \left\{ \sum_{j=2}^{d} \sum_{k=1}^{j-1} \rho_{jk}^2 - \frac{d(d-1)}{2(n-1)} \right\}, \tag{B2}$$

and show that $T$ converges in distribution to the standard normal as $n$ and $d$ grow. The difference is that Mao (2016) uses the exact standard deviation $\sigma_{nd}$, while Leung & Drton (2017) use $d/n$ as an approximation. Leung & Drton (2017) also provide a general theory that applies to other $U$-statistics.

In this simulation, we compare three tests based on Spearman's rho, i.e., the Spearman test in the main paper, the test based on $S$ of Mao (2016), and the test based on $T$ of Leung & Drton (2017).

We apply the three tests on the ten data generating mechanisms described in Section B·2. In additional, we adopt a simulation scheme where data are drawn from independent Cauchy distribution with mean zero and scale one as in Mao (2016) to examine the sizes of the three tests under infinite variance.

Results averaged over $5,000$ simulated data sets are shown in Table 1. The two tests of Mao (2016) and Leung & Drton (2017) have comparable performances across all settings, which agrees with the findings in Mao (2016). We note that the Spearman test achieves higher power against the sparse alternative than the other two tests. This is because our proposed test is based on the maxima while the other two tests are based on averages, and thus our proposed test is more sensitive to the sparse alternatives. We also note that our proposed test can sometimes be conservative, which is a result of the slow convergence rate of the Gumbel distribution. As we will see in Section B·4, this can be addressed by resorting to the simulation-based rejection threshold.

### B·4. *Testing with exact distributions*

In what follows, we provide the empirical sizes and powers of exact tests. We adopt the Gaussian distribution in Section 5·2 in the main paper. We compare the performances of the Spearman test and the Kendall test using theoretical thresholds to the performance of the Spearman and Kendall tests using simulation-based thresholds. We refer to the Spearman test and the Kendall test using simulation-based thresholds as the Spearman exact test and Kendall exact test, respectively.

Results over $5,000$ simulated data sets are given in Table 2. We observe that the sizes of the two exact tests are well controlled, and their powers are higher than the corresponding tests that use the theoretical threshold $q_\alpha$. This reflects the extra gain in power by resorting to the exact tests.

### B·5. *Real data analysis*

We study the empirical performance of competing tests on a real stock market data. We collect the daily closing prices of 452 stocks in the Standard and Poor 500 index from January 1, 2003 to January 1, 2008, available on `finance.yahoo.com`. We study the nearly independent monthly log return data (Xue et al., 2012). All together, the corresponding data matrix has $n = 59$ rows and $d = 452$ columns.

In order to evaluate the control of size for the seven tests, we simulate data sets with independent columns based on the real monthly log return data matrix. We generate each simulated data set by randomly permuting the entries within each column of the data matrix. This permutation preserves the empirical marginal distribution for each of the 452 column variables, i.e. the stock prices, but, within each row, the 452 column variables are mutually independent.

We apply the six competing tests to $1,000$ permuted data sets, and report the resulting $p$-values in Figure 1.

Since the entries within each column have been permuted, the corresponding 452 entries are completely independent and the histograms shall be close to that of the uniform distribution in $[0, 1]$. We find that the histograms of our proposed tests are relatively flat and the proposed tests can effectively control the size. In comparison, the histograms of $p$-values from Zhou (2007) and Mao (2014) are strongly skewed to the left, indicating that the tests tend to falsely reject the null hypothesis. The reason is that Zhou (2007) and Mao (2014) are very sensitive to extreme events as observed in Section 5·2 as well as in Shao & Zhou

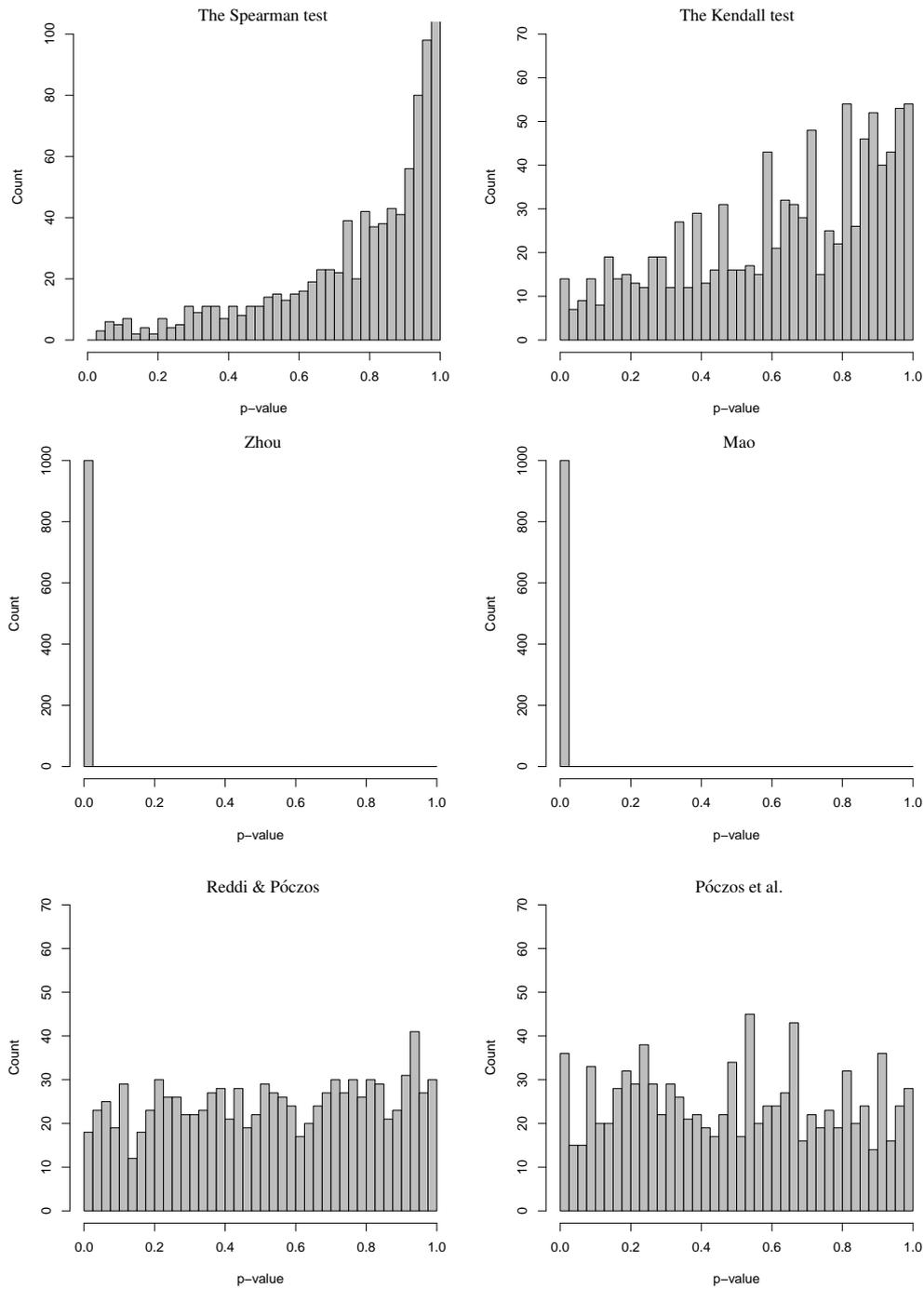

Fig. 1: Histograms of the $p$-values of six competing methods on $1,000$ permuted monthly log return data. The empirical probabilities of the $p$-values less than $0\cdot050$ are $0\cdot003$, $0\cdot021$, $1\cdot000$, $1\cdot000$, $0\cdot041$, and $0\cdot051$ for the Spearman test, the Kendall test, the tests of Zhou (2007), Mao (2014), Reddi & Póczos (2013), and Póczos et al. (2012), respectively.

Table 1: Empirical sizes and powers of the Spearman test, the test of Mao (2016), and the test of Leung & Drton (2017) in percentages

| $n$ | $d$ | Spearman | Leung & Drton (2017) | Mao (2016) | Spearman | Leung & Drton (2017) | Mao (2016) |
|---|---|---|---|---|---|---|---|
| | | | Guassian null distribution | | | Gaussian alternative distribution | |
| 60 | 50 | 2·8 | 5·1 | 5·4 | 91·9 | 30·9 | 31·8 |
| | 200 | 1·8 | 5·1 | 5·2 | 84·3 | 7·4 | 7·5 |
| | 800 | 1·2 | 4·9 | 5·0 | 76·3 | 5·6 | 5·7 |
| 100 | 50 | 3·8 | 4·6 | 4·9 | 97·1 | 59·9 | 60·3 |
| | 200 | 2·5 | 4·6 | 4·7 | 93·7 | 11·6 | 11·8 |
| | 800 | 1·8 | 5·2 | 5·3 | 92·3 | 5·4 | 5·4 |
| | | | Light-tailed null distribution | | | Light-tailed alternative distribution | |
| 60 | 50 | 2·5 | 4·4 | 4·6 | 90·9 | 31·9 | 32·6 |
| | 200 | 1·7 | 4·8 | 5·0 | 84·5 | 6·2 | 6·3 |
| | 800 | 1·1 | 4·8 | 4·9 | 76·0 | 5·3 | 5·4 |
| 100 | 50 | 3·5 | 4·4 | 4·8 | 96·7 | 60·0 | 60·6 |
| | 200 | 2·8 | 5·2 | 5·3 | 94·7 | 10·4 | 10·5 |
| | 800 | 1·8 | 4·7 | 4·8 | 91·7 | 5·9 | 5·9 |
| | | | Heavy-tailed null distribution | | | Heavy-tailed alternative distribution | |
| 60 | 50 | 2·5 | 5·1 | 5·4 | 91·0 | 31·5 | 32·0 |
| | 200 | 1·8 | 5·2 | 5·3 | 84·0 | 7·4 | 7·5 |
| | 800 | 1·1 | 4·2 | 4·3 | 76·0 | 5·3 | 5·4 |
| 100 | 50 | 3·7 | 4·7 | 4·8 | 96·7 | 60·3 | 61·0 |
| | 200 | 3·0 | 4·1 | 4·2 | 94·5 | 11·6 | 11·7 |
| | 800 | 2·1 | 4·7 | 4·8 | 90·9 | 5·2 | 5·3 |
| | | | Multivariate $t$ null distribution | | | Multivariate $t$ alternative distribution | |
| 60 | 50 | 2·8 | 4·4 | 4·6 | 95·2 | 28·7 | 29·5 |
| | 200 | 1·6 | 4·8 | 5·0 | 79·4 | 7·0 | 7·1 |
| | 800 | 1·2 | 5·2 | 5·3 | 40·0 | 5·2 | 5·2 |
| 100 | 50 | 4·1 | 4·7 | 5·1 | 99·7 | 61·0 | 61·6 |
| | 200 | 2·6 | 5·2 | 5·3 | 99·5 | 9·4 | 9·6 |
| | 800 | 1·9 | 4·6 | 4·6 | 98·6 | 5·1 | 5·1 |
| | | | Exponential null distribution | | | Exponential alternative distribution | |
| 60 | 50 | 1·7 | 4·7 | 5·0 | 90·5 | 94·1 | 94·4 |
| | 200 | 0·8 | 5·0 | 5·2 | 83·0 | 100·0 | 100·0 |
| | 800 | 0·2 | 4·3 | 4·4 | 74·7 | 100·0 | 100·0 |
| 100 | 50 | 2·9 | 4·7 | 5·1 | 96·9 | 98·3 | 98·3 |
| | 200 | 1·8 | 5·0 | 5·1 | 94·0 | 100·0 | 100·0 |
| | 800 | 0·7 | 5·3 | 5·4 | 91·4 | 100·0 | 100·0 |
| | | | Cauchy null distribution | | | | |
| 60 | 50 | 1·5 | 4·7 | 4·8 | - | - | - |
| | 200 | 0·6 | 4·7 | 4·8 | - | - | - |
| | 800 | 0·2 | 4·8 | 4·9 | - | - | - |
| 100 | 50 | 3·1 | 5·1 | 5·4 | - | - | - |
| | 200 | 1·7 | 5·1 | 5·1 | - | - | - |
| | 800 | 0·7 | 4·7 | 4·7 | - | - | - |

Results are averaged over 5,000 simulated data sets.

(2014). And here the log return data contain extreme events and are heavy-tailed (Rachev, 2003), which are not eliminated by permutation. Finally, kernel-based tests can control the size, which agrees with our findings in Section 5·2 in the main paper.

## C. TECHNICAL PROOFS
### C·1. *Overview*

In this section, we provide the technical proofs of the theoretical results in the main paper and in Section A of the Supplementary Material. For ease of reading, we defer the technical lemmas to Section C·8.

Table 2: Empirical sizes and powers of simulation-based rejection thresholds in percentages

| $n$ | $d$ | Spearman exact | Kendall exact | Spearman | Kendall | Spearman exact | Kendall exact | Spearman | Kendall |
|---|---|---|---|---|---|---|---|---|---|
| | | \multicolumn{4}{c}{Guassian null distribution} | \multicolumn{4}{c}{Gaussian alternative distribution} | | | | |
| 60 | 50 | 5·6 | 5·4 | 1·8 | 2·9 | 89·9 | 90·7 | 91·9 | 92·8 |
| | 200 | 4·8 | 4·0 | 0·8 | 2·5 | 89·0 | 88·8 | 84·3 | 87·2 |
| | 800 | 4·1 | 4·8 | 0·2 | 1·5 | 97·0 | 97·1 | 97·1 | 97·5 |
| 100 | 50 | 5·9 | 5·8 | 2·8 | 3·7 | 84·5 | 84·4 | 76·3 | 81·8 |
| | 200 | 4·6 | 5·3 | 1·5 | 2·7 | 95·3 | 95·2 | 93·7 | 94·3 |
| | 800 | 5·0 | 4·8 | 0·8 | 2·2 | 94·4 | 94·2 | 92·3 | 93·2 |

The Spearman exact and Kendall exact tests use simulation-based rejection thresholds. Results are averaged over 5,000 simulated data sets.

### C·2. *Proofs of Theorems* 1 *and* 2

In the proof, Lemma C2 plays a key role in calculating the convergence rate of the limiting distribution. We first prove Theorem 1 in the main paper.

*Proof.* To begin with, we focus on the statistic $\psi_{jk} \equiv n^{1/2} V_{jk}/\sigma_V$. In Lemma C2, let $I \equiv \{(j,k) : 1 \le j < k \le d\}$. For $u = (j,k) \in I$, set $B_u = \{(l,m) \in I : (l,m) \ne (j,k), \{l,m\} \cap \{j,k\} \ne \emptyset\}$, $\eta_u = |\psi_{jk}|$, and $A_u = A_{jk} = \{|\psi_{jk}| > t\}$. We can check that $b_3 = 0$ in Lemma C2, and

$$|\mathrm{pr}(n^{1/2} L_n/\sigma_V \le t) - e^{-\lambda_n}| \le b_{1,n} + b_{2,n}, \tag{C1}$$

where we have

$$\lambda_n = \frac{d(d-1)}{2}\mathrm{pr}(A_{12}). \tag{C2}$$

Using Lemma C4, $A_{12}$ is independent of $A_{13}$ and accordingly

$$b_{1,n} \le d^3 \mathrm{pr}(A_{12})^2, \quad b_{2,n} \le d^3 \mathrm{pr}(A_{12} A_{13}) = d^3 \mathrm{pr}(A_{13})^2.$$

Here using Lemma C5, when $t = o(n^{1/6})$, we have

$$\mathrm{pr}(A_{12}) = \mathrm{pr}(|\psi_{12}| > t) = 2\{1 - \Phi(t)\}\{1 + o(1)\}. \tag{C3}$$

Accordingly, for $i = 1, 2$, using the Gaussian tail bound $\mathrm{pr}\{N_1(0,1) > t\} \le e^{-t^2/2}/\{(2\pi)^{1/2} t\}$, we have

$$b_{i,n} \le \frac{2}{\pi t^2} d^3 \exp(-t^2) \Rightarrow b_{1,n} + b_{2,n} \le \frac{4}{\pi t^2} d^3 \exp(-t^2)\{1 + o(1)\}. \tag{C4}$$

We then let

$$t = (4\log d - \log\log d + y)^{1/2} \asymp (4\log d)^{1/2}, \tag{C5}$$

and directly plug the above $t$ into (C1). Because $\log d = o(n^{1/3})$, (C3) holds and it follows that

$$b_{1n} + b_{2n} \le \frac{4}{\pi(4\log d - \log\log d + y)} d^3 \exp(-4\log d + \log\log d) = o\!\left(\frac{1}{d}\right). \tag{C6}$$

On the other hand, using the Gaussian tail bounds in an unpublished technical report by Duembgen (available on `arXiv.org` with identifier 1012.2063), we have for any $t > 0$,

$$\frac{1}{t + 1/t}(2\pi)^{-1/2}\exp\!\left(-\frac{t^2}{2}\right) \le 1 - \Phi(t) \le \frac{1}{t}(2\pi)^{-1/2}\exp\!\left(-\frac{t^2}{2}\right). \tag{C7}$$

Accordingly, as $d$ grows, we see that $t$ diverges to infinity in (C5). We have, as $t$ grows,

$$1/t - 1/(t + 1/t) = 1/\{t(t^2 + 1)\} \asymp 1/t^3.$$

It yields that

$$1 - \Phi(t) = \frac{1}{(2\pi)^{1/2} t} \exp\left(-\frac{t^2}{2}\right) [1 + O\{(\log d)^{-3/2}\}]. \tag{C8}$$

Combining (C2), (C3), and (C8) implies

$$\lambda_n = d^2 \{1 - \Phi(t)\}\{1 + o(1)\} = \frac{d^2}{(8\pi \log d)^{1/2}} \exp\left(-\frac{4 \log d - \log \log d + y}{2}\right)\{1 + o(1)\}$$

$$= (8\pi)^{-1/2} \exp\left(-\frac{y}{2}\right)\{1 + o(1)\}. \tag{C9}$$

Plugging the above equation to (C1) yields

$$\left|\text{pr}\left(\frac{nL_n^2}{\sigma_V^2} - 4\log d + \log \log d \le y\right) - \exp\left\{-(8\pi)^{-1/2} \exp\left(-\frac{y}{2}\right)\right\}\right|$$

$$\le |\text{pr}(n^{1/2} L_n/\sigma_V \le t) - \exp(-\lambda_n)| + |\exp(-\lambda_n) - \exp\{-(8\pi)^{-1/2} \exp(-y/2)\}| = o_y(1), \tag{C10}$$

which completes the proof. □

The proof of Theorem 2 is very similar to the proof of Theorem 1. One only needs to replace (C22) with (C23) when applying Lemma C5. The proof is thus omitted.

### C·3. *Proofs of Theorems* 3 *and* 4

The proofs are based on several concentration inequalities developed in Section C·8. We prove Theorem 3 first.

*Proof.* The test statistic $nL_n^2/\sigma_V^2$ is scale and location invariant. Hence, without loss of generality, we assume that $\sum_{i=1}^n c_{ni} = 0$ in this proof. Using (4), we have $E_{H_0}(V_{jk}) = 0$ and

$$\widehat{V}_{jk} = \frac{V_{jk}}{\sigma_V} = \frac{V_{jk}\{1 + o(1)\}}{A_1}.$$

Let $\Delta$ be a Lipschitz constant of both $g(\cdot)$ introduced in (1) and $f(\cdot)$ introduced in (2) in the main paper. Using Lemma C6, it follows that, for sufficiently large $n$ and some scalar $c(A_1, A_2, \Delta)$ only depending on $A_1, A_2$, and $\Delta$, for any $t > 0$,

$$\text{pr}(|\widehat{V}_{jk} - V_{jk}| > t) \le 2\exp\{-nt^2/c(A_1, A_2, \Delta)\}.$$

We then have

$$\text{pr}\left(\max_{j<k} |\widehat{V}_{jk} - V_{jk}| > t\right) \le d^2 \exp\{-nt^2/c(A_1, A_2, \Delta)\},$$

which implies that, with probability at least $1 - d^{-1}$,

$$\max_{j,k} |\widehat{V}_{jk} - V_{jk}| \le \left\{\frac{3c(A_1, A_2, \Delta) \log d}{n}\right\}^{1/2}.$$

Therefore, we have, for $n$ large enough, there exists a large enough constant $C$ such that

$$nL_n^2/\sigma_V^2 = n\max_{j<k} \widehat{V}_{jk}^2 \ge n\left(\max_{j,k} |V_{jk}| - \max_{j,k} |\widehat{V}_{jk} - V_{jk}|\right)^2 \ge \{C - 3^{1/2} c(A_1, A_2, \Delta)^{1/2}\}^2 \log d.$$

Accordingly, by choosing $C > 2 + 3^{1/2} c(A_1, A_2, \Delta)^{1/2}$, we have with probability no smaller than $1 - d^{-1}$,

$$nL_n^2/\sigma_V^2 > (4 + \epsilon) \log d,$$

for some small constant $\epsilon$. Accordingly, for any given $q_\alpha$, with probability tending to 1,

$$nL_n^2/\sigma_V^2 > 4\log d - \log \log d - q_\alpha.$$

This completes the proof. □

We then prove Theorem 4 in the main paper.

*Proof.* The proof is similar to that of Theorem 3. Because the test statistic $n\widetilde{L}_n^2/\sigma_U^2$ is scale and location invariant, without loss of generality, we assume $E_{H_0}\{h(X_1,\ldots,X_m)\} = 0$. Then it is immediately clear that $E_{H_0}(U_{jk}) = 0$. Moreover, by a standard argument of $U$-statistics (see, e.g., Serfling (2002)), we have

$$\begin{aligned}
n\mathrm{var}_{H_0}(U_{jk}) &= \widetilde{\sigma}_U^2\{1+o(1)\} \\
&= m^2\mathrm{var}_{H_0}[E_{H_0}\{h(X_{1,\{1,2\}},\ldots,X_{m,\{1,2\}}) \mid X_{1,\{1,2\}}\}]\{1+o(1)\} \\
&= A_4\{1+o(1)\},
\end{aligned}$$

where $\widetilde{\sigma}_U^2$ is defined in (15) in the main paper. Then using Lemma C7, we have for large $n$ and some scalar $c(A_3, A_4, m)$ only depending on $A_3, A_4$ and $m$, for any $t > 0$

$$\mathrm{pr}\big(|\widehat{U}_{jk} - U_{jk}| > t\big) \leq 2\exp\{-nt^2/c(A_3, A_4, m)\}.$$

The rest is a line-by-line follow of Theorem 3's proof. □

C·4. *Proof of Theorem* 5

*Proof.* Consider the Gaussian setting and a simple alternative set of parameters

$$\mathcal{F}(\rho) = \{M : M = I_d + \rho e_1 e_j^{\mathrm{T}} + \rho e_j e_1^{\mathrm{T}}, e_k = (\underbrace{0,\ldots,0}_{k-1}, 1, 0, \ldots, 0), 1 \leq k \leq d, j = 2,\ldots,d\}.$$

Let $\mu_\rho$ be the uniform measure on $\mathcal{F}(\rho)$ and $\rho = c_0(\log d/n)^{1/2}$ for some small enough constant $c_0 < 1$. Let $\mathrm{pr}_\Sigma$ denote the probability measure of $N_d(0,\Sigma)$ and $\mathrm{pr}_{\mu_\rho} = \int \mathrm{pr}_\Sigma \mathrm{d}\mu_\rho(\Sigma)$. Let $\mathrm{pr}_0$ denote the probability measure of $N_d(0, I_d)$. Note that, for any set $A$, we have

$$\sup_{\Sigma \in \mathcal{F}(\rho)} \mathrm{pr}_\Sigma(A^C) \geq \mathrm{pr}_{\mu_\rho}(A^C), \quad 1 = \mathrm{pr}_{\mu_\rho}(A^C) + \mathrm{pr}_{\mu_\rho}(A),$$

and

$$\mathrm{pr}_{\mu_\rho}(A) \leq \mathrm{pr}_0(A) + |\mathrm{pr}_{\mu_\rho}(A) - \mathrm{pr}_0(A)|.$$

Letting $A \equiv \{T_\alpha = 1\}$, the above equations yield

$$\inf_{T_\alpha \in \mathcal{T}_\alpha} \sup_{\Sigma \in \mathcal{F}(\rho)} \mathrm{pr}_\Sigma(T_\alpha = 0) \geq 1 - \alpha - \sup_{A:\mathrm{pr}_0(A)\leq\alpha} |\mathrm{pr}_{\mu_\rho}(A) - \mathrm{pr}_0(A)| \geq 1 - \alpha - \frac{1}{2}\|\mathrm{pr}_{\mu_\rho} - \mathrm{pr}_0\|_{TV},$$

where $\|\cdot\|_{TV}$ denotes the total variation norm. Setting $L_{\mu_\rho}(y) \equiv \mathrm{dpr}_{\mu_\rho}(y)/\mathrm{dpr}_0(y)$, and by Jensen's inequality, we have

$$\|\mathrm{pr}_{\mu_\rho} - \mathrm{pr}_0\|_{TV} = \int |L_{\mu_\rho}(y) - 1|\mathrm{dpr}_0(y) = E_{\mathrm{pr}_0}|L_{\mu_\rho}(Y) - 1| \leq [E_{\mathrm{pr}_0}\{L_{\mu_\rho}^2(Y)\} - 1]^{1/2}.$$

Therefore, as long as $E_{\mathrm{pr}_0}\{L_{\mu_\rho}^2(Y)\} = 1 + o(1)$, we have

$$\inf_{T_\alpha \in \mathcal{T}_\alpha} \sup_{\Sigma \in \mathcal{F}(\rho)} \mathrm{pr}_\Sigma(T_\alpha = 0) \geq 1 - \alpha - o(1) > 0. \tag{C11}$$

We then prove that $E_{\mathrm{pr}_0}\{L_{\mu_\rho}^2(Y)\} = 1 + o(1)$. By construction, we have

$$L_{\mu_\rho} = \frac{1}{d-1}\sum_{\Sigma\in\mathcal{F}(\rho)}\left[\prod_{i=1}^n \frac{1}{|\Sigma|^{1/2}}\exp\left\{-\frac{1}{2}Z_{i,\cdot}^{\mathrm{T}}(\Omega - I_d)Z_{i,\cdot}\right\}\right],$$

where $\Omega \equiv \Sigma^{-1}$ and $Z_{1,\cdot}, \ldots, Z_{n,\cdot}$ are $d$-dimensional vectors to be specified later. We have

$$E_{\mathrm{pr}_0}\{L^2_{\mu_\rho}(Y)\} = \frac{1}{(d-1)^2} \sum_{\Sigma_1, \Sigma_2 \in \mathcal{F}(\rho)} E\left[\prod_{i=1}^n \frac{1}{|\Sigma_1|^{1/2}} \frac{1}{|\Sigma_2|^{1/2}} \exp\left\{-\frac{1}{2} Z_{i,\cdot}^{\mathrm{T}} (\Omega_1 + \Omega_2 - 2I_d) Z_{i,\cdot}\right\}\right],$$

where $\Omega_i \equiv \Sigma_i^{-1}$ for $i = 1, 2$ and $\{Z_{i,\cdot}, 1 \leq i \leq n\}$ are independent and identically distributed as $N_d(0, I_d)$. We write

$$A = \frac{\rho}{1-\rho^2} \begin{pmatrix} 2\rho & -1 & -1 \\ -1 & \rho & 0 \\ -1 & 0 & \rho \end{pmatrix}, \quad B = \frac{2\rho}{1-\rho^2} \begin{pmatrix} \rho & -1 \\ -1 & \rho \end{pmatrix}.$$

It is easy to derive that

$$E_{\mathrm{pr}_0}(L^2_{\mu_\rho}) = \underbrace{\frac{d-2}{d-1} \prod_{i=1}^n \left[\frac{1}{1-\rho^2} E\left\{\exp\left(-\frac{1}{2} Z_{i,\{1,2,3\}}^{\mathrm{T}} A Z_{i,\{1,2,3\}}\right)\right\}\right]}_{E_1}$$

$$+ \underbrace{\frac{1}{d-1} \prod_{i=1}^n \left[\frac{1}{1-\rho^2} E\left\{\exp\left(-\frac{1}{2} Z_{i,\{1,2\}}^{\mathrm{T}} B Z_{i,\{1,2\}}\right)\right\}\right]}_{E_2},$$

where $E_1$ represents the set of $(\Sigma_1, \Sigma_2)$ with $\Sigma_1 \neq \Sigma_2$, and $E_2$ represents the set of $(\Sigma_1, \Sigma_2)$ with $\Sigma_1 = \Sigma_2$. By standard argument in moment generating functions of the Gaussian quadratic form, we have

$$E_1 = \frac{d-2}{d-1} \frac{1}{(1-\rho^2)^n} \left[\{1 + \lambda_1(A)\}\{1 + \lambda_2(A)\}\{1 + \lambda_3(A)\}\right]^{-n/2},$$

where $\lambda_i(A)$ is the $i$th eigenvalue of $A$. Moreover, we have $\{1 + \lambda_1(A)\}\{1 + \lambda_2(A)\}\{1 + \lambda_3(A)\} = |A + I_d| = (1-\rho^2)^{-2}$. When $d$ grows with $n$, we know that

$$E_1 = \frac{1}{(1-\rho^2)^n}(1-\rho^2)^n\{1 + o(1)\} = 1 + o(1). \tag{C12}$$

For $E_2$, it is easy to calculate that $\lambda_1(B) = 2\rho/(1-\rho)$ and $\lambda_2(B) = -2\rho/(1+\rho)$. Similar to the calculation of $E_1$, we have $E_2 = (d-1)^{-1}(1-\rho^2)^{-n}$. Recalling that $\rho = c_0(\log d/n)^{1/2}$ and $\log d/n = o(1)$, we have

$$E_2 = (d-1)^{-1}(1 - c_0^2 \log d/n)^{-n} = (d-1)^{-1} \exp(c_0^2 \log d)\{1 + o(1)\} = o(1) \tag{C13}$$

as long as $c_0 < 1$. Combining (C12) and (C13) yields (C11). Lastly, because the Pearson's covariance matrix $\Sigma \in \mathcal{F}(\rho)$ implies that the Pearson's correlation matrix $R \in \mathcal{F}(\rho)$, we have $\{X : \Sigma \in \mathcal{F}(\rho)\} \subset \{X : R \in \mathcal{F}(\rho)\}$ and thus

$$\inf_{T_\alpha \in \mathcal{T}_\alpha} \sup_{R \in \mathcal{F}(\rho)} \mathrm{pr}_\Sigma(T_\alpha = 0) \geq \inf_{T_\alpha \in \mathcal{T}_\alpha} \sup_{\Sigma \in \mathcal{F}(\rho)} \mathrm{pr}_\Sigma(T_\alpha = 0) \geq 1 - \alpha - o(1) > 0.$$

This completes the proof. □

### C·5. *Proofs of Theorems A3 and A5*

We first prove Theorem A3.

*Proof.* By checking the proof of Theorem 4 in Cai & Jiang (2011), we only need to verify the following three statements to show that Theorem A3 holds.

**S1.** Suppose that $Z \equiv (Z_1, Z_2, Z_3, Z_4)^{\mathrm{T}} \sim N_4(0, \Sigma_1)$ with

$$\Sigma_1 \equiv \begin{bmatrix} 1 & 0 & r & 0 \\ 0 & 1 & 0 & 0 \\ r & 0 & 1 & 0 \\ 0 & 0 & 0 & 1 \end{bmatrix}, \quad |r| \leq 1.$$

Let $Z_{1,\cdot}, \ldots, Z_{n,\cdot} \in \mathcal{R}^4$, with $Z_{i,\cdot} = (Z_{i,1}, \ldots, Z_{i,4})^{\mathrm{T}}$, be $n$ independent observations of $Z$. Further set $t_n \equiv \{(4\log d - \log\log d + y)/n\}^{1/2}$ for some fixed $y \in \mathcal{R}$, as $n$ grows, and $\log d = o(n^{1/3})$. We have

$$\sup_{|r| \leq 1} \mathrm{pr}(3\tau_{12}/2 > t_n, 3\tau_{34}/2 > t_n) = O(d^{-4}),$$

where for $(j, k) \in \{(1, 2), (3, 4)\}$,

$$\tau_{jk} \equiv \frac{2}{n(n-1)} \sum_{1 \leq i < i' \leq n} \mathrm{sign}(Z_{i,j} - Z_{i',j})\mathrm{sign}(Z_{i,k} - Z_{i',k}).$$

**S2** Suppose that $Z \equiv (Z_1, Z_2, Z_3, Z_4)^{\mathrm{T}} \sim N_4(0, \Sigma_2)$ with

$$\Sigma_2 \equiv \begin{bmatrix} 1 & 0 & r_1 & 0 \\ 0 & 1 & r_2 & 0 \\ r_1 & r_2 & 1 & 0 \\ 0 & 0 & 0 & 1 \end{bmatrix}, \quad |r_1| \leq 1, |r_2| \leq 1.$$

Let $Z_{1,\cdot}, \ldots, Z_{n,\cdot} \in \mathcal{R}^4$, with $Z_{i,\cdot} = (Z_{i,1}, \ldots, Z_{i,4})^{\mathrm{T}}$, be $n$ independent observations of $Z$. Then set $t_n \equiv \{(4\log d - \log\log d + y)/n\}^{1/2}$ for some fixed $y \in \mathcal{R}$, $n$ and $d$ grow, and $\log d = o(n^{1/3})$. We have

$$\sup_{|r_1| \leq 1, |r_2| \leq 1} \mathrm{pr}(3\tau_{12}/2 > t_n, 3\tau_{34}/2 > t_n) = O(d^{-4}).$$

**S3** Suppose that $Z \equiv (Z_1, Z_2, Z_3, Z_4)^{\mathrm{T}} \sim N_4(0, \Sigma_3)$ with

$$\Sigma_3 \equiv \begin{bmatrix} 1 & 0 & r_1 & 0 \\ 0 & 1 & 0 & r_2 \\ r_1 & 0 & 1 & 0 \\ 0 & r_2 & 0 & 1 \end{bmatrix}, \quad |r_1| \leq 1, |r_2| \leq 1.$$

Let $Z_{1,\cdot}, \ldots, Z_{n,\cdot} \in \mathcal{R}^4$, with $Z_{i,\cdot} = (Z_{i,1}, \ldots, Z_{i,4})^{\mathrm{T}}$, be $n$ independent replicates of $Z$. Then setting $t_n \equiv \{(4\log d - \log\log d + y)/n\}^{1/2}$ for some fixed $y \in \mathcal{R}$, as $n$ and $d$ grow, and $\log d = o(n^{1/3})$. Then we have, for any fixed $\delta \in (0, 1)$, there exists $\epsilon_0 = \epsilon(\delta) > 0$ such that

$$\sup_{|r_1|, |r_2| \leq 1-\delta} \mathrm{pr}(3\tau_{12}/2 > t_n, 3\tau_{34}/2 > t_n) = O(d^{-2-\epsilon_0}).$$

For showing **S1, S2,** and **S3** hold, consider the general setting where $Z \equiv (Z_1, Z_2, Z_3, Z_4)^{\mathrm{T}} \sim N_4(0, \Sigma)$ and $\Sigma$ has diagonals all equal one. The Kendall's tau correlation coefficient is a $U$-statistic with degree two and the kernel function bounded by one. By exploiting the Hájek's projection (Hájek et al., 1999), with a little abuse of notation, we can write

$$3\tau_{jk}/2 = \frac{2}{n} \sum_{i=1}^{n} E(3\tau_{jk}/2 \mid Z_{i,\{j,k\}}) + E_{jk} = \frac{1}{n} \sum_{i=1}^{n} \underbrace{E(3\tau_{jk} \mid Z_{i,\{j,k\}})}_{\Psi_{i,jk}} + E_{jk}, \quad (\text{C14})$$

where $\Psi_{1,jk}, \Psi_{2,jk}, \ldots, \Psi_{n,jk}$ are $n$ independent random variables, and $E_{jk}$ is a degenerate $U$-statistic. Moreover, both $\Psi_{i,jk}$ and $E_{jk}$ are bounded. Using (C14) and the Slutsky's argument, we can further write

$$\begin{aligned}&\operatorname{pr}(3\tau_{12}/2 > t_n, 3\tau_{34}/2 > t_n)\\ =&\operatorname{pr}\Big(\frac{1}{n}\sum_{i=1}^{n}\Psi_{i,12} + E_{12} > t_n, \frac{1}{n}\sum_{i=1}^{n}\Psi_{i,34} + E_{34} > t_n\Big)\\ \leq&\operatorname{pr}\Big(\frac{1}{n}\sum_{i=1}^{n}\Psi_{i,12} > t_n - \epsilon_1, \frac{1}{n}\sum_{i=1}^{n}\Psi_{i,34} > t_n - \epsilon_1\Big) + \operatorname{pr}(E_{12} > \epsilon_1) + \operatorname{pr}(E_{34} > \epsilon_1)\\ =&\operatorname{pr}\Big\{n^{-1/2}\sum_{i=1}^{n}\Psi_{i,12} > n^{1/2}(t_n - \epsilon_1), n^{-1/2}\sum_{i=1}^{n}\Psi_{i,34} > n^{1/2}(t_n - \epsilon_1)\Big\}\\ &+ \operatorname{pr}(E_{12} > \epsilon_1) + \operatorname{pr}(E_{34} > \epsilon_1)\end{aligned}$$

where $\epsilon_1$ is a constant to be specified later. Because $|\Psi_{i,jk} n^{-1/2}| \leq 3n^{-1/2}$ for $(j,k) \in \{(1,2),(3,4)\}$, using Theorem 1 in Zaïtsev (1987), we have

$$\begin{aligned}&\operatorname{pr}\Big\{n^{-1/2}\sum_{i=1}^{n}\Psi_{i,12} > n^{1/2}(t_n - \epsilon_1), n^{-1/2}\sum_{i=1}^{n}\Psi_{i,34} > n^{1/2}(t_n - \epsilon_1)\Big\}\\ \leq&\operatorname{pr}\Big\{Y_1 \geq n^{1/2}(t_n - \epsilon_1 - \epsilon_2), Y_2 \geq n^{1/2}(t_n - \epsilon_1 - \epsilon_2)\Big\} + c_1 \exp\Big(-n\epsilon_2/c_2\Big),\end{aligned} \qquad (C15)$$

where $c_1$ and $c_2$ are two positive constants and $(Y_1, Y_2)^{\mathrm{T}}$ is bivariate Gaussian with mean zero and covariance matrix

$$\Sigma_Y = \operatorname{cov}\Big\{\Big(n^{1/2}\sum_{i=1}^{n}\Psi_{i,12}, n^{1/2}\sum_{i=1}^{n}\Psi_{i,34}\Big)^{\mathrm{T}}\Big\}.$$

We then determine what $\Sigma_Y$ is. Recall that under **S1, S2,** or **S3**, $Z_j, Z_k$ are independent for $(j,k) \in \{(1,2),(3,4)\}$. We can write

$$\Psi_{i,jk} = E(3\tau_{jk} \mid Z_{i,\{j,k\}}) = 3E\{\operatorname{sign}(Z_{i,j} - \widetilde{Z}_j)\operatorname{sign}(Z_{i,k} - \widetilde{Z}_k) \mid Z_{i,j}, Z_{i,k}\},$$

where $(\widetilde{Z}_j, \widetilde{Z}_k)^{\mathrm{T}}$ is an independent copy of $(Z_{i,j}, Z_{i,k})^{\mathrm{T}}$. Because $\widetilde{Z}_j$ is independent of $\widetilde{Z}_k$, we can write

$$\begin{aligned}&3E\{\operatorname{sign}(Z_{i,j} - \widetilde{Z}_j)\operatorname{sign}(Z_{i,k} - \widetilde{Z}_k) \mid Z_{i,j}, Z_{i,k}\}\\ =&3E\{\operatorname{sign}(Z_{i,j} - \widetilde{Z}_j) \mid Z_{i,j}\}E\{\operatorname{sign}(Z_{i,k} - \widetilde{Z}_k) \mid Z_{i,k}\}\\ =&3\{\operatorname{pr}(\widetilde{Z}_j > Z_{i,j} \mid Z_{i,j}) - \operatorname{pr}(\widetilde{Z}_j < Z_{i,j} \mid Z_{i,j})\}\{\operatorname{pr}(\widetilde{Z}_k > Z_{i,k} \mid Z_{i,k}) - \operatorname{pr}(\widetilde{Z}_k < Z_{i,k} \mid Z_{i,k})\}.\end{aligned} \qquad (C16)$$

Using the property of the Gaussian distribution, (C16) yields

$$\Psi_{i,jk} = 3\{1 - 2\Phi(Z_{i,j})\}\{1 - 2\Phi(Z_{i,k})\}, \qquad (C17)$$

where $\Phi(\cdot)$ is the distribution function of the standard Gaussian. Using the result in Example 2 in the main paper, we know

$$n\operatorname{var}(\tau_{jk}) = \frac{2(2n+5)}{9(n-1)} = \frac{4}{9} + o(1).$$

Combining it with Lemma A in Page 183 in Serfling (2002) yields that

$$n\operatorname{var}(3\tau_{jk}) = 4 + o(1) = 4\operatorname{var}(\Psi_{1,jk}) + O(n^{-1}).$$

Because $\operatorname{var}(\Psi_{i,jk})$ is a constant irrelevant to $n$, we have $\operatorname{var}(\Psi_{i,jk}) = 1$ for $i = 1, \ldots, n$ and $(j,k) \in \{(1,2),(3,4)\}$. This yields

$$[\Sigma_Y]_{11} = [\Sigma_Y]_{22} = \operatorname{var}(\Psi_{1,12}) = 1.$$

In the end, we determine the value of $[\Sigma_Y]_{12}$. It is immediately clear that

$$[\Sigma_Y]_{12} = \mathrm{cov}\Big(n^{-1/2}\sum_{i=1}^n \Psi_{i,12}, n^{-1/2}\sum_{i=1}^n \Psi_{i,34}\Big) = \mathrm{cov}(\Psi_{1,12}, \Psi_{1,34}).$$

Using (C17), we can further write

$$\mathrm{cov}(\Psi_{1,12}, \Psi_{1,34}) = 9E\{1 - 2\Phi(Z_1)\}\{1 - 2\Phi(Z_2)\}\{1 - 2\Phi(Z_3)\}\{1 - 2\Phi(Z_4)\}. \tag{C18}$$

Using (C18), we are now ready to prove that statements **S1, S2,** and **S3** hold. Recall that $(Y_1, Y_2)^\mathrm{T} \sim N_2(0, I_2)$. (C14) yields

$$\mathrm{pr}(3\tau_{12}/2 > t_n, 3\tau_{34}/2 > t_n) \leq \mathrm{pr}\{Y_1 \geq n^{1/2}(t_n - \epsilon_1 - \epsilon_2), Y_2 \geq n^{1/2}(t_n - \epsilon_1 - \epsilon_2)\} \\ + c_1 \exp(-n\epsilon_2/c_2) + \mathrm{pr}(E_{12} > \epsilon_1) + \mathrm{pr}(E_{34} > \epsilon_1).$$

Both $E_{12}$ and $E_{34}$ are degenerate $U$-statistics with kernel function bounded. From Proposition 2.3 in Arcones & Gine (1993), we know that there exist constants $c_3, c_4$ such that

$$\mathrm{pr}(E_{12} > \epsilon_1) \leq c_3 \exp(-c_4 n\epsilon_1), \quad \mathrm{pr}(E_{34} > \epsilon_1) \leq c_3 \exp(-c_4 n\epsilon_1).$$

Recalling that $t_n = \{(4\log d - \log\log d + y)/n\}^{1/2} \asymp (4\log d/n)^{1/2}$ and $\log d = o(n^{1/3})$, we can pick $\epsilon_1, \epsilon_2$ small enough such that $\epsilon_1, \epsilon_2 \asymp n^{-2/3}$. In this way, we have for any constant $c > 0$, there exists a scalar $C$ depending on $c$ such that, for $n$ large enough,

$$\exp(-cn\epsilon_i) \leq \exp(-Cn^{1/3}) = o(d^{-4}), \quad i = 1, 2,$$

and $\epsilon_1 = o(t_n), \epsilon_2 = o(t_n)$.

For **S1** and **S2**, we know that $Z_4$ is independent of $Z_1, Z_2, Z_3$, and accordingly

$$\mathrm{cov}(\Psi_{1,12}, \Psi_{1,34}) = 9E\{1 - 2\Phi(Z_1)\}\{1 - 2\Phi(Z_2)\}\{1 - 2\Phi(Z_3)\}\{1 - 2\Phi(Z_4)\} \\ = 9E\{1 - 2\Phi(Z_1)\}\{1 - 2\Phi(Z_2)\}\{1 - 2\Phi(Z_3)\}E\{1 - 2\Phi(Z_4)\} = 0.$$

Therefore, we have $(Y_1, Y_2)^\mathrm{T} \sim N_2(0, I_2)$ and accordingly

$$\mathrm{pr}(3\tau_{12}/2 > t_n, 3\tau_{34}/2 > t_n) \leq [\mathrm{pr}\{Y_1 \geq n^{1/2}(t_n - \epsilon_1 - \epsilon_2)\}]^2 + c_1 \exp(-n\epsilon_2/c_2) \\ + \mathrm{pr}(E_{12} > \epsilon_1) + \mathrm{pr}(E_{34} > \epsilon_1) \\ = (\mathrm{pr}[Y_1 \geq n^{1/2}t_n\{1 + o(1)\}])^2 + o(d^{-4}) = o(d^{-4}),$$

where we use the Gaussian tail bound that for any $t > 0$,

$$\{\mathrm{pr}(Y_1 > t)\}^2 \leq \frac{2}{\pi t^2}\exp(-t^2).$$

For proving **S3**, we need one more lemma, which shows that $[\Sigma_Y]_{12}$ is upper bounded by a constant strictly less than 1 when all off-diagonal values in $\Sigma_3$ are upper bounded by $r < 1$.

LEMMA C1. *Suppose that* $(Z_1, Z_2, Z_3, Z_4)^\mathrm{T} \sim N_4(0, \Sigma_{\mathrm{full}})$ *with*

$$\Sigma_{\mathrm{full}} = \begin{bmatrix} 1 & a_1 & a_2 & a_3 \\ a_1 & 1 & a_4 & a_5 \\ a_2 & a_4 & 1 & a_6 \\ a_3 & a_5 & a_6 & 1 \end{bmatrix}.$$

*If* $|a_1|, |a_2|, \ldots, |a_6| \leq r < 1$, *then we have*

$$\sup_{|a_1|, |a_2|, \ldots, |a_6| \leq r} |\mathrm{corr}[\{\Phi(Z_1) - 1/2\}\{\Phi(Z_2) - 1/2\}, \{\Phi(Z_3) - 1/2\}\{\Phi(Z_4) - 1/2\}]| = C_r < 1.$$

*Here $C_r \leq 1$ only depends on $r$. Moreover, we have $C_r = 1$ only when $r = 1$ and $\{a_1, a_2, \ldots, a_6\}$ attain the boundary that $|a_j| = 1$ for some $j \in \{1, \ldots, 6\}$.*

*Proof.* First, we show that $C_r = 1$ only when $r = 1$ and $\{a_1, a_2, \ldots, a_6\}$ attain the boundary. When $C_r = 1$, we have

$$\{\Phi(Z_1) - 1/2\}\{\Phi(Z_2) - 1/2\} = a\{\Phi(Z_3) - 1/2\}\{\Phi(Z_4) - 1/2\}$$

for some constant $a$. This implies that

$$Z_1 = \Phi^{-1}\left[\frac{a\{\Phi(Z_3) - 1/2\}\{\Phi(Z_4) - 1/2\}}{\Phi(Z_2) - 1/2} + 1/2\right].$$

We have $Z_1 \sim N_1(0,1)$ if and only if $a\{\Phi(Z_3) - 1/2\}\{\Phi(Z_4) - 1/2\}/\{\Phi(Z_2) - 1/2\} \sim \mathrm{Unif}(-1/2, 1/2)$. Here $\mathrm{Unif}(-1/2, 1/2)$ represents the random variable uniformly distributed in the interval $[-1/2, 1/2]$. Because when $Z_2 \neq \pm Z_3$ and $Z_2 \neq \pm Z_4$, there is always possibility such that $Z_2$ is very close to zero and both $Z_3$ and $Z_4$ are away from zero, such that $a\{\Phi(Z_3) - 1/2\}\{\Phi(Z_4) - 1/2\}/\{\Phi(Z_2) - 1/2\}$ is very close to $\infty$ and outside of $[-1/2, 1/2]$. Accordingly, $Z_2$ must be equal to either $\pm Z_3$ or $\pm Z_4$. Or equivalently, $\{a_1, a_2, \ldots, a_6\}$ attain the boundary $r = 1$. This completes the proof of the first part.

Secondly, it is obvious that there is a one-to-one mapping between $r$ and

$$C_r \equiv \sup_{|a_1|, |a_2|, \ldots, |a_6| \leq r} |\mathrm{corr}[\{\Phi(Z_1) - 1/2\}\{\Phi(Z_2) - 1/2\}, \{\Phi(Z_3) - 1/2\}\{\Phi(Z_4) - 1/2\}]|.$$

Accordingly, as long as $r < 1$, $C_r < 1$ only depends on $r$. □

Using Lemma C1, we can continue to prove **S3** holds. Recall that now $(Y_1, Y_2)^\mathrm{T} \sim N_2(0, \Sigma_Y)$, where Lemma C1 shows $\sup_{|r_1|, |r_2| \leq 1-\delta} |[\Sigma_Y]_{12}| \leq C_r < 1$. Thus, we have

$$\mathrm{pr}(Y_1 \geq t, Y_2 \geq t) = \mathrm{pr}\{\min(Y_1, Y_2) \geq t\}.$$

Denoting $\rho \equiv [\Sigma_Y]_{12}$, using Equation (8) in Nadarajah & Kotz (2008), we have

$$E[\exp\{t \min(Y_1, Y_2)\}] = \exp\left(\frac{t^2}{2}\right)\Phi\left\{\frac{-t(1-\rho)}{(2-2\rho)^{1/2}}\right\}.$$

Using the Chernoff's bounding method, we immediately have

$$\sup_{|r_1|, |r_2| \leq 1-\delta} \mathrm{pr}(Y_1 \geq t, Y_2 \geq t) \leq \sup_{|r_1|, |r_2| \leq 1-\delta} \inf_{\lambda > 0} \frac{E[\exp\{\lambda \min(Y_1, Y_2)\}]}{e^{\lambda t}}$$

$$\leq \sup_{|r_1|, |r_2| \leq 1-\delta} \inf_{\lambda > 0} e^{\lambda^2/2 - \lambda t}\Phi\left\{\frac{-\lambda(1-\rho)}{(2-2\rho)^{1/2}}\right\}$$

$$= \inf_{\lambda > 0} e^{\lambda^2/2 - \lambda t}\Phi\left[-\lambda\{(1-C_r)/2\}^{1/2}\right].$$

Picking $\lambda = t$, the above equation yields

$$\sup_{|r_1|, |r_2| \leq 1-\delta} \mathrm{pr}(Y_1 \geq t, Y_2 \geq t) \leq e^{-t^2/2}\Phi\left[-t\{(1-C_r)/2\}^{1/2}\right].$$

Setting $t = n^{1/2} t_n \{1 + o(1)\}$, then there exists a constant $C$ such that

$$\sup_{|r_1|, |r_2| \leq 1-\delta} \mathrm{pr}(Y_1 \geq t, Y_2 \geq t) \leq Cd^{-2}(\log d)^{1/2} \mathrm{pr}[Y_1 > \{(1-C_r)/2\}^{1/2} t]$$

$$\leq Cd^{-2}(\log d)^{1/2} O(d^{-M}), \tag{C19}$$

where $M > 0$ is a constant only depending on $C_r$. Thus, the statement **S3** holds.

All in all, we have **S1, S2,** and **S3** all hold. This completes the proof. □

We then proceed to prove Theorem A5.

*Proof.* We strictly follow the proof of Theorem 5 in the main paper and adopt the same notation system. In particular, we consider the following alternative set of parameters:

$$\mathcal{F}_m(\rho) = \{\Sigma^0 = I_d + \rho e_1 e_j^{\mathrm{T}} + \rho e_j e_1^{\mathrm{T}}, \text{ for } j \in \{m+1, m+2, \ldots, d\}\}.$$

Then the whole proof in Theorem 5 applies here with the only exception that $E_2 = (d-m)^{-1}(1-\rho^2)^{-n}$. However, because $d - m \asymp d$, we have $E_2 = (d-m)^{-1}(1-\rho^2)^{-n} \asymp d^{-1}(1-\rho^2)^{-n}$. Taking $\rho = c_0'(\log d/n)^{1/2}$, we still have

$$E_2 \asymp \frac{1}{d}(1 - c_0'^2 \log d/n)^n = d^{-1} \exp(c_0'^2 \log d)\{1 + o(1)\} = o(1).$$

This completes the proof. □

### C·6. *Proofs of Theorems A6, A7, and A8*

The proof of Theorem A6 is very similar to that of Theorem 1 and is accordingly omitted. In the following we give the proof of Theorem A7.

*Proof.* Assume that the first entry across $X_{1,\cdot}, \ldots, X_{n,\cdot}$ is heterogeneity. It is obvious that $\mathrm{sign}(X_{i,1} - X_{i',1})$ is invariant to $\beta_0$ and $\sigma^2$ given $\beta_1/\sigma$. Therefore, without loss of generality, we assume $\beta_0 = 0$ and $\sigma^2 = 1$. Moreover, without loss of generality, we can assume $\beta_1 \in (0, M)$, otherwise we can always replace $\beta_1$ with $\min(|\beta_1|, M)$. We have

$$\begin{aligned}
E\{\mathrm{sign}(X_{i',1} - X_{i,1})\} &= \mathrm{pr}(X_{i',1} - X_{i,1} > 0) - \mathrm{pr}(X_{i',1} - X_{i,1} < 0) \\
&= \mathrm{pr}\{Z_{i',1} - Z_{i,1} > -\beta_1(i' - i)/n\} - \mathrm{pr}\{Z_{i',1} - Z_{i,1} < -\beta_1(i' - i)/n\} \\
&= \mathrm{pr}\{Z_{i',1} - Z_{i,1} < \beta_1(i' - i)/n\} - \mathrm{pr}\{Z_{i',1} - Z_{i,1} < -\beta_1(i' - i)/n\},
\end{aligned}$$

where $Z_{k,1} = \{X_{k,1} - E(X_{k,1})\}/\mathrm{var}(X_{k,1})$ is the standardized version of $Z_{k,1}$ for $k = 1, \ldots, n$. Then, (A3) yields that the density function of $Z_{i',1} - Z_{i,1}$ is

$$\begin{aligned}
\{p_{i'1} * (-p_{i1})\}(z) &= \int_{-\infty}^{\infty} p_{i'1}(z+y)p_{i1}(y)\mathrm{d}y \geq D_4 \int_{-M}^{M} p_{i'1}(z+y)\mathrm{d}y \\
&\geq D_4 \int_{\max\{-M+z, -M\}}^{\min\{M+z, M\}} p_{i'1}(y)\mathrm{d}y \geq D_4^2(2M - |z|), \quad (|z| \leq 2M).
\end{aligned}$$

This further implies

$$\begin{aligned}
E(h_1) &= \frac{2}{n(n-1)} \sum_{i<i'} E\{\mathrm{sign}(X_{i',1} - X_{i,1})\} = \frac{2}{n(n-1)} \sum_{i<i'} [F_p\{\beta_1(i'-i)/n\} - F_p\{-\beta_1(i'-i)/n\}] \\
&\geq \frac{2D_4^2 M \beta}{n^2(n-1)} \sum_{i<i'} (i'-i) = \frac{D_4^2 M \beta}{3} \frac{n(n^2-1)}{n^2(n-1)} \geq \frac{2D_4^2 M \beta}{3},
\end{aligned}$$

where $F_p(\cdot)$ is the distribution function of $Z_{1,1} - Z_{2,1}$. On the other hand, by the McDiarmid's inequality (McDiarmid, 1989), for any $j \in \{1, \ldots, d\}$,

$$\mathrm{pr}\{|h_j - E(h_j)| > t\} \leq 2\exp(-nt^2/2).$$

The rest is similar to the proof of Theorem 3 in the main paper. □

We then proceed to prove Theorem A8.

*Proof.* We focus on a simple Gaussian model where $X_{1,\cdot}, \ldots, X_{n,\cdot}$ are independent and normally distributed, with covariance matrix $I_d$. Accordingly, by virtue of the normal distribution, we can write $(X_{1,j}, \ldots, X_{n,j})^{\mathrm{T}} \sim N_n(\mu_{j,\cdot}, I_n)$ for $j \in \{1, \ldots, d\}$. Here $\mu_{j,\cdot} \in \mathcal{R}^n$ is the mean vector. We then consider the following simple alternative set of parameters:

$$\mathcal{H}(\beta) = \Big\{\mu = \{\mu_1, \ldots, \mu_d\} : \mu_{i,\cdot} = \{0, \beta/n, 2\beta/n, \ldots, (n-1)\beta/n\}^{\mathrm{T}} \text{ for some } i, \text{ the rests are all zero}\Big\}.$$

Let $\mu_\beta$ be the uniform measure on $\mathcal{H}(\beta)$ and $\beta = c_0''(\log d/n)^{1/2}$ for some small enough constant $c_0'' < 3^{1/2}$. Let $\mathrm{pr}_\mu$ be the probability measure on $N_n(\mu_{1,\cdot}, I_n) \otimes \cdots \otimes N_n(\mu_{n,\cdot}, I_n)$. In particular, let $\mathrm{pr}_0$ be the probability measure on $N_n(0, I_n) \otimes \cdots \otimes N_n(0, I_n)$. Let $\mathrm{pr}_{\mu_\beta} \equiv \int \mathrm{pr}_\mu \mathrm{d}\mu_\beta(\mu)$ be the measure based on $\mathcal{H}(\beta)$. Similar to the proof of Theorem 5, to prove Theorem A8, it suffices to show that

$$E_{\mathrm{pr}_0}\{L_{\mu_\beta}^2(Y)\} = 1 + o(1),$$

where $L_{\mu_\beta}(y) \equiv \mathrm{dpr}_{\mu_\beta}(y)/\mathrm{dpr}_0(y)$. By construction, we can write

$$L_{\mu_\beta}(y) = \frac{1}{d} \sum_{\mu \in \mathcal{H}(\beta)} \Big\{\prod_{i=1}^d \exp\Big(Z_{i,\cdot}^{\mathrm{T}} \mu_{i,\cdot} - \|\mu_{i,\cdot}\|_2^2/2\Big)\Big\}.$$

Accordingly, the above equation yields that

$$E_{\mathrm{pr}_0}\{L_{\mu_\beta}^2(Y)\} = \frac{1}{d^2} \sum_{\mu^1, \mu^2 \in \mathcal{H}(\beta)} E\Big\{\prod_{i=1}^d \exp(Z_{i,\cdot}^{\mathrm{T}} \mu_{i,\cdot}^1 + Z_{i,\cdot}^{\mathrm{T}} \mu_{i,\cdot}^2 - \|\mu_{i,\cdot}^1\|_2^2/2 - \|\mu_{i,\cdot}^2\|_2^2/2)\Big\},$$

where $Z_{1,\cdot}, \ldots, Z_{d,\cdot} \sim N_n(0, I_n)$ and $\mu^k = \{\mu_{1,\cdot}^k, \ldots, \mu_{d,\cdot}^k\}$ for $k \in \{1, 2\}$. We can then continue to write

$$E_{\mathrm{pr}_0} L_{\mu_\beta}^2 = \underbrace{\frac{1}{d^2} \sum_{\mu^1 \neq \mu^2} E\Big\{\prod_{i=1}^d \exp(Z_{i,\cdot}^{\mathrm{T}} \mu_{i,\cdot}^1 + Z_{i,\cdot}^{\mathrm{T}} \mu_{i,\cdot}^2 - \|\mu_{i,\cdot}^1\|_2^2/2 - \|\mu_{i,\cdot}^2\|_2^2/2)\Big\}}_{H_1} +$$

$$\underbrace{\frac{1}{d^2} \sum_{\mu^1 = \mu^2} E\Big\{\prod_{i=1}^d \exp(Z_{i,\cdot}^{\mathrm{T}} \mu_{i,\cdot}^1 + Z_{i,\cdot}^{\mathrm{T}} \mu_{i,\cdot}^2 - \|\mu_{i,\cdot}^1\|_2^2/2 - \|\mu_{i,\cdot}^2\|_2^2/2)\Big\}}_{H_2}. \qquad (C20)$$

Let $\mu^* \equiv \{0, \beta/n, \ldots, (n-1)\beta/n\}^{\mathrm{T}}$. For the first term in (C20), we have

$$H_1 = \frac{d-1}{d} E\{\exp(Z_{1,\cdot}^{\mathrm{T}} \mu^* - \|\mu^*\|_2^2/2)\} E\{\exp(Z_{2,\cdot}^{\mathrm{T}} \mu^* - \|\mu^*\|_2^2/2)\} = 1 + o(1).$$

For the second term in (C20), we have, when $c_0'' \leq \sqrt{3}$,

$$H_2 = d^{-1} E\left\{\exp(2Z_{1,\cdot}^{\mathrm{T}} \mu^* - \|\mu^*\|_2^2)\right\} = d^{-1} \exp(\|\mu^*\|_2^2) = d^{-1} \exp\{(1 - n^{-1})(2n-1)\beta^2/6\}$$
$$= d^{-1} \exp(n\beta^2/3)\{1 + o(1)\} = \exp\{-\log d + (c_0'')^2 \log d/3\}\{1 + o(1)\} = o(1).$$

This completes the proof. □

### C·7. *Proof of Theorem A1*

*Proof.* We focus on simple linear rank statistics, as the extension to rank-type $U$-statistics is straightforward. Following the proof of Theorem 1 in the main paper and using Lemma C5, we can replace (C3) with

$$\mathrm{pr}(A_{12}) = \mathrm{pr}(|\psi_{12}| > t) = 2\{1 - \Phi(t)\}[1 + O\{(\log d)^{3/2} n^{-1/2} + (\log d)^{1/2} n^{-1/6}\}].$$

Furthermore, (C8) implies that

$$\lambda_n = d^2\{1 - \Phi(t)\}[1 + O\{(\log d)^{3/2} n^{-1/2} + (\log d)^{1/2} n^{-1/6}\}][1 + O\{(\log d)^{-3/2}\}]$$
$$= (8\pi)^{-1/2} \exp\Big(-\frac{y}{2}\Big)[1 + O\{(\log d)^{3/2} n^{-1/2} + (\log d)^{1/2} n^{-1/6} + (\log d)^{-3/2}\}].$$

Accordingly, we can separately bound the first and second terms in (C10), yielding that

$$|\mathrm{pr}(L_n \leq t) - \exp(-\lambda_n)| = o(d^{-1})$$

and

$$|\exp(-\lambda_n) - \exp\{-\exp(-y/2)(8\pi)^{-1/2}\}| = O\{(\log d)^{3/2}n^{-1/2} + (\log d)^{1/2}n^{-1/6} + (\log d)^{-3/2}\}.$$

Here we use the fact that, when $x$ approaches zero, $\exp(x) - 1 \asymp x$. This completes the proof. □

### C·8. *Auxiliary lemmas*

The following seven lemmas play crucial roles in our theory.

LEMMA C2 (ARRATIA ET AL. (1989)). *Let $I$ be an index set and $\{B_\alpha, \alpha \in I\}$ be a set of subsets of $I$; that is, $B_\alpha \subset I$ for each $\alpha \in I$. Let also $\{\eta_\alpha, \alpha \in I\}$ be random variables. For a given $t \in \mathcal{R}$, set $\lambda = \sum_{\alpha \in I} \mathrm{pr}(\eta_\alpha > t)$. Then*

$$\left|\mathrm{pr}\big(\max_{\alpha \in I} \eta_\alpha \leq t\big) - e^{-\lambda}\right| \leq \min(1, \lambda^{-1})(b_1 + b_2 + b_3),$$

*where*

$$b_1 \equiv \sum_{\alpha \in I} \sum_{\beta \in B_\alpha} \mathrm{pr}(\eta_\alpha > t)\mathrm{pr}(\eta_\beta > t), \quad b_2 \equiv \sum_{\alpha \in I} \sum_{\beta \neq \alpha, \beta \in B_\alpha} \mathrm{pr}(\eta_\alpha > t, \eta_\beta > t),$$

$$b_3 \equiv \sum_{\alpha \in I} E|\mathrm{pr}\{\eta_\alpha > t \mid \sigma(\eta_\beta, \beta \notin B_\alpha)\} - \mathrm{pr}(\eta_\alpha > t)|,$$

*where $\sigma(\eta_\beta, \beta \notin B_\alpha)$ is the $\sigma$-algebra generated by $\{\eta_\beta, \beta \notin B_\alpha\}$. In particular, if $\eta_\alpha$ is independent of $\{\eta_\beta, \beta \notin B_\alpha\}$ for each $\alpha$, then $b_3 = 0$.*

LEMMA C3. *Suppose that $X, Y$ are two independent continuous random variables. Let $X_1, \ldots, X_n$ and $Y_1, \ldots, Y_n$ be independent observations of $X$ and $Y$. Let $\{Q_i^X, i = 1, \ldots, n\}$ and $\{Q_i^Y, i = 1, \ldots, n\}$ be the rank of $X_i$ and $Y_i$ in the samples $\{X_i\}_{i=1}^n$ and $\{Y_i\}_{i=1}^n$. Let $\{R_{ni}\}_{i=1}^n$ represent the relative ranks:*

$$R_{ni} = Q_{i'}^Y \quad \text{subject to} \quad Q_{i'}^X = i.$$

*We then have $\{R_{n1}, \ldots, R_{nn}\}$ are uniformly distributed in all permutations of $\{1, \ldots, n\}$ with*

$$\mathrm{pr}(R_{n1} = i_1, \ldots, R_{nn} = i_n) = \frac{1}{n!}, \tag{C21}$$

*for any permeation $\{i_1, \ldots, i_n\}$ of $\{1, \ldots, n\}$. Here $n!$ represents the factorial of $n$.*

*Proof.* Using the fact that $\{X_i\}_{i=1}^n$ are independent of $\{Y_i\}_{i=1}^n$, for any permutation $\{i_1, \ldots, i_n\}$ of $\{1, \ldots, n\}$ and any $a_1, \ldots, a_n \in \mathcal{R}$, we have

$$\mathrm{pr}(X_{i_1} < X_{i_2} < \cdots < X_{i_n} \mid Y_1 = a_1, \ldots, Y_n = a_n) = \mathrm{pr}(X_{i_1} < X_{i_2} < \cdots < X_{i_n}).$$

Therefore, the relative ranks' joint distribution is identical to the distribution of $\{Q_i^X, i = 1, \ldots, n\}$. The latter's distribution is known to be jointly distributed in the form of (C21). □

LEMMA C4. *Let $\{S_{jk}, 1 \leq j < k \leq d\}$ be functions of relative ranks $\{R_{ni}^{jk}, i = 1, \ldots, n\}$ with the same mapping function from $\{R_{ni}^{jk}, i = 1, \ldots, n\}$ for any $j, k$. Then, under the null hypothesis $H_0$, $S_{u_1 j}$ is identically and pairwise independently distributed to $S_{u_2 k}$ for any non-identical $(u_1, j)$ and $(u_2, k)$.*

*Proof.* Using Lemma C3, the distribution of the relative ranks does not change as long as the independence assumption holds. We then have $\{S_{jk}, 1 \leq j < k \leq d\}$ are all identically distributed. It is obvious that, under $H_0$, $S_{u_1 j}, S_{u_2 k}$ are independent when there is no overlap between $(u_1, j)$ and $(u_2, k)$. In the rest we show that $S_{u_1 j}, S_{u_2 k}$ are independent when there is one overlap between $(u_1, j)$ and $(u_2, k)$.

We consider the case $u_1 = u_2 \neq j \neq k$ and the proofs of all the other settings are similar. We prove $S_{uj}$ is independent of $S_{uk}$ with $u = u_1 = u_2 \in \{1, \ldots, d\}$. It is sufficient to show that for any two bounded and measurable functions $g(x)$ and $h(x)$, we have

$$E\{g(S_{uj})h(S_{uk})\} = E\{g(S_{uj})\}E\{h(S_{uk})\}.$$

Given $\{X_{1,u}, X_{2,u}, \ldots, X_{n,u}\}$, $S_{uj}$ and $S_{uk}$ are independent. We have

$$E\{g(S_{uj})h(S_{uk})\} = E(E[g(S_{uj})h(S_{uk}) \mid \{X_{1,u}, X_{2,u}, \ldots, X_{n,u}\}])$$
$$= E(E[g(S_{uj}) \mid \{X_{1,u}, X_{2,u}, \ldots, X_{n,u}\}]E[h(S_{uk}) \mid \{X_{1,u}, X_{2,u}, \ldots, X_{n,u}\}]).$$

Next we show that, given $\{X_{1,u}, X_{2,u}, \ldots, X_{n,u}\}$, the conditional distributions of $S_{uj}$ and $g(S_{uj})$ are irrelevant to $\{X_{1,u}, X_{2,u}, \ldots, X_{n,u}\}$. This follows by applying Lemma C3. A detailed proof can be found in Pages 477–479 in Kendall & Stuart (1961). Using this argument, we then have

$$E[g(S_{uj}) \mid \{X_{1,u}, X_{2,u}, \ldots, X_{n,u}\}] = E[g(S_{uj}) \mid \{X'_{1,u}, X'_{2,u}, \ldots, X'_{n,u}\}],$$

for any sequence $\{X'_{1,u}, X'_{2,u}, \ldots, X'_{n,u}\}$ randomly drawn from $X_u$. This implies

$$E[g(S_{uj}) \mid \{X_{1,u}, X_{2,u}, \ldots, X_{n,u}\}] = E\{g(S_{uj})\}.$$

Similarly, we have

$$E[g(S_{uk}) \mid \{X_{1,u}, X_{2,u}, \ldots, X_{n,u}\}] = E\{g(S_{uk})\}.$$

This shows that $\{S_{jk}, 1 \leq j < k \leq d\}$ are pairwise independent. □

LEMMA C5. *Suppose that the boundedness assumption in Theorem 2 hold. We then have, in a region $x \in (0, o(n^{1/6}))$,*

$$\text{pr}\left[\frac{U_{jk} - E(U_{jk})}{\{\text{var}(U_{jk})\}^{1/2}} > x\right] = \{1 - \Phi(x)\}\left\{1 + O\left(\frac{1 + x^3}{n^{1/2}}\right)\right\}. \tag{C22}$$

*Suppose that the regularity conditions in Theorem 1 hold. Under the null hypothesis $H_0$ holds, we have in the region $x \in (0, O(n^{1/6-\epsilon}))$ for some $\epsilon > 0$,*

$$\text{pr}\left[\frac{V_{jk} - E_{H_0}(V_{jk})}{\{\text{var}_{H_0}(V_{jk})\}^{1/2}} > x\right] = \{1 - \Phi(x)\}\left\{1 + O\left(\frac{1 + x^3}{n^{1/2}} + \frac{x}{n^{1/6}}\right)\right\}. \tag{C23}$$

*And we can replace the rate in the right-hand side of* (C23) *with $1 + o(1)$ when we have $x \in (0, o(n^{1/6}))$.*

*Proof.* For the moderate deviation properties of the $U$-statistics, the general results for them of unbounded kernel functions can be found in Malevich & Abdalimov (1979) and Vandemaele (1983). Borovskikh & Weber (2003) give the result for $U$-statistics of bounded kernels with symmetric kernels. However, using a similar argument as in Eichelsbacher (1998) and Hoeffding (1948), the results can be generalized to the multivariate data and asymmetric kernel cases.

When we do not specify the rate of convergence on the right hand side of (C23), the proof of the moderate deviation for simple linear rank statistics is in Kallenberg (1982). For explicitly characterizing the rate, we simply follow Kallenberg (1982). Below we adopt some notation used in Kallenberg (1982). Consider the data with $n$ independent samples $X_1, \ldots, X_n$ drawn from $X \in \mathcal{R}$. Let $F(\cdot)$ be the distribution function of $X$. Let $R_{n1}, \ldots, R_{nn}$ be the ranks of $X_1, \ldots, X_n$. Let $S_n = \sum_{i=1}^n c_{ni} g\{R_{ni}/(n+1)\}$ be the simple linear rank statistic of interest and $V_n = \sum_{i=1}^n c_{ni} g\{F(X_i)\}$ be an intermediate one. It is obvious that $S_n$ is identically distributed to $V_{jk}$ under the null hypothesis.

Let $\mu_n$ and $\tau_n$ be the mean and standard deviation of $S_n$. Without loss of generality, we assume $\mu_n = 0$. Equations (2.1) and (2.2) in Kallenberg (1982) imply

$$\text{pr}(S_n > x\tau_n) \geq \text{pr}\{V_n > (x + n^{-1/6}\tau_n)\} - \text{pr}(|S_n - V_n| > n^{-1/6}\tau_n)$$
$$\text{and } \text{pr}(S_n > x\tau_n) \leq \text{pr}\{V_n > (x - n^{-1/6}\tau_n)\} + \text{pr}(|S_n - V_n| > n^{-1/6}\tau_n). \tag{C24}$$

On one hand, using the lemma in Page 406 in Kallenberg (1982), we have

$$\text{pr}(|S_n - V_n| > n^{-1/6}\tau_n)\{1 - \Phi(x)\}^{-1} \leq (1/2)^{\delta n^{1/3}}\{1 - \Phi(n^{1/6-\epsilon})\}^{-1}$$
$$\leq \exp\{-(\delta n^{1/3})\log 2 + n^{1/3-2\epsilon}/2\}O(n^{1/6-\epsilon}). \tag{C25}$$

On the other hand, $V_n$ is the sum of independent bounded random variables. Therefore, we can use the classic result on the moderate deviation of sums of independence variables (check, for example, Chapter 8 in Petrov (1975)). It implies that for any $y_n$,

$$\text{pr}(V_n > y_n \tau_n) = \{1 - \Phi(y_n)\}\Big\{1 + O\Big(\frac{1 + y_n^3}{n^{1/2}}\Big)\Big\}. \tag{C26}$$

We let $|y_n - x| \leq n^{-1/6}$, which implies that $1 + y_n^3 \asymp 1 + x^3$. Then, standard arguments on Gaussian tail probabilities give us

$$\{1 - \Phi(y_n)\}/\{1 - \Phi(x)\} = 1 + O(n^{-1/6}x). \tag{C27}$$

Plugging (C25), (C26), and (C27) into (C24), we have

$$\text{pr}(S_n > x\tau_n) = \{1 - \Phi(x)\}\Big\{1 + O\Big(\frac{1 + y_n^3}{n^{1/2}} + \frac{x}{n^{1/6}}\Big)\Big\}.$$

This completes the proof. □

LEMMA C6 (CONCENTRATION INEQUALITY FOR SIMPLE LINEAR RANK STATISTICS). *Assume the setting and notation in Lemma* 3. *Consider the simple linear rank statistic*

$$V \equiv \sum_{i=1}^{n} c_{ni} g\Big(\frac{R_{ni}}{n+1}\Big) = \frac{1}{n}\sum_{i=1}^{n} f\Big(\frac{Q_i^X}{n+1}\Big) g\Big(\frac{Q_i^Y}{n+1}\Big),$$

*where $f(\cdot)$ and $g(\cdot)$ are Lipschitz functions with Lipschitz constant $\Delta < \infty$ and $\max\{|f(0)|, |g(0)|\} \leq A_2$. We have, for any $t > 0$,*

$$\text{pr}(|V - EV| > t) \leq 2\exp(-Cnt^2),$$

*for some scalar $C$ only depending on $\Delta$ and $A_2$.*

*Proof.* The proof is an application of the McDiarmid's inequality (McDiarmid, 1989). In the samples $\{(X_i, Y_i), i = 1, \ldots, n\}$, consider replacing $(X_1, Y_1)$ with $(X_1', Y_1')$ and fix all the others. Then the ranks of $\{Q_i^X, i = 1, \ldots, n\}$ and $\{Q_i^Y, i = 1, \ldots, n\}$ are changed to $\{\widetilde{Q}_i^X, i = 1, \ldots, n\}$ and $\{\widetilde{Q}_i^Y, i = 1, \ldots, n\}$. By the alignment assumption, we have

$$\Big|\sum_{i=1}^{n} c_{ni} g\Big(\frac{R_{ni}}{n+1}\Big) - \sum_{i=1}^{n} c_{ni} g\Big(\frac{\widetilde{R}_{ni}}{n+1}\Big)\Big| = \frac{1}{n}\Big|\sum_{i=1}^{n} f\Big(\frac{Q_i^X}{n+1}\Big) g\Big(\frac{Q_i^Y}{n+1}\Big) - \sum_{i=1}^{n} f\Big(\frac{\widetilde{Q}_i^X}{n+1}\Big) g\Big(\frac{\widetilde{Q}_i^Y}{n+1}\Big)\Big|.$$

Because $\max_{1 \leq i \leq n} |f\{i/(n+1)\}| \leq A_2 + \Delta$ and $\max_{1 \leq i \leq n} |g\{i/(n+1)\}| \leq A_2 + \Delta$, it yields

$$\frac{1}{n}\Big|\sum_{i=1}^{n} f\Big(\frac{Q_i^X}{n+1}\Big) g\Big(\frac{Q_i^Y}{n+1}\Big) - \sum_{i=1}^{n} f\Big(\frac{\widetilde{Q}_i^X}{n+1}\Big) g\Big(\frac{\widetilde{Q}_i^Y}{n+1}\Big)\Big|$$
$$\leq \frac{A_2 + \Delta}{n}\Big\{\sum_{i=1}^{n}\Big|f\Big(\frac{Q_i^X}{n+1}\Big) - f\Big(\frac{\widetilde{Q}_i^X}{n+1}\Big)\Big| + \sum_{i=1}^{n}\Big|g\Big(\frac{Q_i^Y}{n+1}\Big) - g\Big(\frac{\widetilde{Q}_i^Y}{n+1}\Big)\Big|\Big\}.$$

Here the inequality follows from the fact that for any two sequences $\{(x_1^1, y_1^1), \ldots, (x_n^1, y_n^1)\}$ and $\{(x_1^2, y_1^2), \ldots, (x_n^2, y_n^2)\}$,

$$\Big|\sum_{i=1}^{n} x_i^1 y_i^1 - \sum_{i=1}^{n} x_i^2 y_i^2\Big| \leq \sum_{i=1}^{n} |x_i^1(y_i^1 - y_i^2)| + \sum_{i=1}^{n} |y_i^2(x_i^1 - x_i^2)|$$
$$\leq \max_{1 \leq i \leq n} |x_i^1| \sum_{i=1}^{n} |y_i^1 - y_i^2| + \max_{1 \leq i \leq n} |y_i^2| \sum_{i=1}^{n} |x_i^1 - x_i^2|.$$

Using the fact that both $f(\cdot)$ and $g(\cdot)$ are Lipschitz, we can further write

$$\frac{A_2 + \Delta}{n}\Big\{\sum_{i=1}^n \Big|f\Big(\frac{Q_i^X}{n+1}\Big) - f\Big(\frac{\widetilde{Q}_i^X}{n+1}\Big)\Big| + \sum_{i=1}^n \Big|g\Big(\frac{Q_i^Y}{n+1}\Big) - g\Big(\frac{\widetilde{Q}_i^Y}{n+1}\Big)\Big|\Big\}$$

$$\leq \frac{\Delta(A_2 + \Delta)}{n(n+1)}\Big(\sum_{i=1}^n |Q_i^X - \widetilde{Q}_i^X| + \sum_{i=1}^n |Q_i - \widetilde{Q}_i^Y|\Big).$$

Because only one position in $\{X_1, \ldots, X_n\}$ and $\{Y_1, \ldots, Y_n\}$ is changing, we have

$$\sum_{i=1}^n |Q_i^X - \widetilde{Q}_i^X| \leq 2(n-1) \text{ and } \sum_{i=1}^n |Q_i^Y - \widetilde{Q}_i^Y| \leq 2(n-1).$$

This further implies that

$$\Big|\sum_{i=1}^n c_{ni}g\Big(\frac{R_{ni}}{n+1}\Big) - \sum_{i=1}^n c_{ni}g\Big(\frac{\widetilde{R}_{ni}}{n+1}\Big)\Big| \leq \frac{4(A_2+\Delta)\Delta(n-1)}{n(n+1)} \asymp \frac{1}{n}.$$

Then, by using the McDiarmid's inequality, we have the desired concentration inequality. $\square$

LEMMA C7 (CONCENTRATION INEQUALITY FOR $U$-STATISTICS). *Suppose that $U$ is a $U$-statistic with degree $m$ and bounded kernel $|h(\cdot)| \leq M$. We then have, for any $t > 0$,*

$$\mathrm{pr}(|U - EU| > t) \leq 2\exp\{-nt^2/(2mM^2)\}.$$

*Proof.* This concentration inequality follows from calculating the moment generating function of the $U$-statistics and using the Hoeffding's decoupling trick. Check Hoeffding (1963) for the detailed proof. $\square$

LEMMA C8. *Under the Gaussian model with the Pearson's correlation matrix R, we have the following four equations hold:*

$$E(\rho_{jk}) = \frac{6}{\pi}\arcsin(R_{jk}/2) + O(1/n), \quad E(\tau_{jk}) = \frac{2}{\pi}\arcsin(R_{jk}),$$

$$E(\widehat{\rho}_{jk}) = \frac{6}{\pi}\arcsin(R_{jk}/2) + O(1/n), \text{ and } E(\widehat{\tau}_{jk}) = \frac{4}{\pi}\arcsin(R_{jk}/2) + O(1/n).$$

*Proof.* The relationship between Spearman's rho, Kendall's tau, and Pearson's correlation coefficients under the Gaussian model can be found in Kruskal (1958). Noticing that $\widehat{\rho}_{jk}$ and $\widehat{\tau}_{jk}$ are asymptotically equivalent to $\rho_{jk}$ and $2\rho_{jk}/3$, we have the other two equations. $\square$